\documentclass[11pt]{article}
\usepackage{amsmath,amssymb,graphicx,float,epsf}
\def\RR{\hbox{I\kern-.2em\hbox{R}}}
\newcommand{\qed}{\hbox to 0pt{}\hfill$\rlap{$\sqcap$}\sqcup$ \vspace{3mm}}
\numberwithin{equation}{section}
\restylefloat{figure}
\usepackage{graphicx} 
\usepackage{amsmath}
\usepackage{amsthm}
\usepackage{amsfonts}
\usepackage[caption = false]{subfig}
\usepackage{wrapfig}
\usepackage{enumerate}
\usepackage{enumitem}
\usepackage{array}
\usepackage{caption}
\usepackage[margin=1cm]{caption}
\usepackage{authblk}

\usepackage[utf8]{inputenc}

\newtheorem{Lemma}{Lemma}
\newtheorem{Th}{Theorem}

\newtheorem{Ex}{Example}

\usepackage{graphicx}
\usepackage{tikz}
\usetikzlibrary{shapes.geometric, arrows}
\tikzstyle{rect} = [draw, rectangle, fill=blue!20, text width=6em, text centered, minimum height=2em]
\tikzstyle{elli} = [draw, ellipse, fill=red!20, minimum height=2em]
\tikzstyle{circ} = [draw, circle, fill=white!20, minimum width=8pt, inner sep=5pt]
\tikzstyle{diam} = [draw, diamond, fill=white!20, text width=6em, text badly centered, inner sep=0pt]
\tikzstyle{line} = [draw, -latex']
\usepackage{a4wide} 
\usepackage{authblk} 
\usepackage[superscript,biblabel]{cite}
\usepackage{varioref}
\usepackage{hyperref}
\hypersetup{colorlinks=true, linkcolor=violet, citecolor=cyan, urlcolor=teal} 

\graphicspath{ {Images2/} }
\date{}
\begin{document}
	\title{Dynamics of Reaction-Diffusion-Advection System and its Impact on River Ecology in the Presence of Spatial Heterogeneity I}
	
	
	\author[1]{\small Md. Kamrujjaman\thanks{Corresponding author: Md Kamrujjaman, E-mail: kamrujjaman@du.ac.bd}}
	\author[1,2]{\small Mayesha Sharmim Tisha\thanks{ Email:  tisha.math@aiub.edu}} 

	\affil[1]{\footnotesize Department of Mathematics, University of Dhaka, Dhaka 1000, Bangladesh}
	\affil[2]{\footnotesize Department of Mathematics, American International University-Bangladesh (AIUB), Kuratoli, Khilkhet, Dhaka 1229, Bangladesh}	

\maketitle  	
\vspace{-1.0cm}
\noindent\rule{6.25in}{0.02in}\\
 {\bf Abstract.}
In this study, a spatially distributed reaction-diffusion-advection (RDA) model with harvesting is investigated to signify the outcome of a competition between two competing species in a heterogeneous environment. 
The study builds upon the concept presented in literature \cite{tisha2}, applying it to river ecology in the context of harvesting activities. We assume that despite of having distinct advection and diffusion rates, two species are competing for the same food supply. This paper's main objective is to study, using theoretical and numerical analysis, the global asymptotic stability and coexistence steady state based on different and unequal rates of diffusion and advection. We establish the result for existence, uniqueness and positivity of the solution. The local stability of two semi trivial steady states is demonstrated. Also, we examine the non-existence of coexistence steady state with the help of some non-trivial presumptions. Finally, we combine the local stability with the non-existence of coexistence to demonstrate the global stability using monotone dynamical systems.
\noindent\rule{6.25in}{0.02in}\\

	 
\clearpage	 
\section{Introduction}
River habitats can be identified by their main unidirectional flow. This flow greatly promotes the spread of several species, including stream insects, algae, and invertebrates. According to M\"{u}ller (1982), the ``drift paradox" \cite{drift} is the question of how a population can endure in a river in spite of washout caused by flow. It refers to the mystery of how populations withstand washout and continue to establish themselves in such advective environment. Recent modelling papers have addressed this problem (Lutscher et al., 2005; Pachepsky et al., 2005; Speirs and Gurney, 2001).

The growth and expansion of a biological population are shaped by a complex interplay of environmental factors and population density. The spatial configuration and characteristics of the environment play a crucial role in guiding the movement patterns of individuals within the population. These movement patterns, whether they are random or influenced by external environmental cues, directly impact the population's overall dynamics and distribution. Individual movement can be broadly categorized into two types: undirected movement, where individuals move randomly without a specific target, and directed movement, where movement is influenced by environmental gradients, resource availability, or other external factors. Understanding these movement behaviors is essential for studying population dynamics and predicting patterns of population spread in varying ecological contexts.	Many species, from gut-dwelling bacteria to river-dwelling flora and fauna, reside in habitats with primarily unidirectional flow, which produces a significant bias in the way that organisms like algae, crustaceans, and stream insects spread. These aquatic species, which reside in streams and rivers and are constantly vulnerable to downstream advection due to water flow, are possibly the most notable example.
	
The reaction-diffusion-advection system is a mathematical framework that can model how organisms move and spread in response to environmental factors like food availability, predation, and habitat conditions \cite{gan22,lee24}. By incorporating spatial heterogeneity-such as variations in river depth, flow velocity, and nutrient levels-this model helps explain the distribution and migration patterns of river animals, allowing for better management of their habitats.
	
	
The structure and flow of a river, determined by factors such as the gradient of the riverbed and the velocity of the current, create distinct zones within river ecosystems. In forested areas, trees contribute the majority of the food supply for streams, while in other regions, algae become the dominant food source. Species in rivers and lakes adapt to specific locations where food availability, light conditions, predation, competition, water depth, and current velocity shape their survival \cite{yan22}.
	
Through floods, landslides, the building of dams and channels, water extraction, or changes in land use patterns within the watershed, human activity has the potential to modify these natural circumstances.
Understanding the effects of heterogeneity on persistence and extinction, invasion potential and invasibility, as well as competition of such populations, is therefore essential for the sustainable management and maintenance of riverine ecosystems \cite{he3}. Drift paradox \cite{drift} refers to the mystery of how populations withstand washout and continue to establish themselves in such advective environment.
	
One of the most crucial factors in population dynamics from an environmental standpoint is the impact of harvesting. Harvesting includes activities like hunting, fishing, or other techniques for lowering population size and density. 
	
The consequences of harvesting, which is a major cause of species extinction, are the most important factor in population ecology when it comes to managing population. It is critical to understand the harvesting threshold in order to make the best use of limited natural resources while also maintaining the ecological system's balance. In the case of any particular ecological resources, population extinction brought on by harvesting should be avoided if the main objective is to achieve sustainable control in population management \cite{prey-taxis,iz2022,cosner14}. 
The primary objective of this study is to investigate the impact of harvesting on a reaction-diffusion-advection model that describes the dynamics of a two-species population in a riverine habitat. This work builds upon and expands the methodologies and concepts presented in our recent publication \cite{tisha2}. In the referenced study, we developed a mathematical model to represent population dynamics within a flowing medium, incorporating key processes such as diffusion, advection, and reaction terms that capture intra- and interspecies interactions. By synthesizing these foundational approaches, the present work delves deeper into how harvesting activities influence the ecological balance and spatial-temporal distribution of the two species. The model formulation in this study maintains the fundamental structure outlined in \cite{tisha2}, with modifications to include harvesting terms, allowing for a comprehensive analysis of its effects on population stability, persistence, and distribution in the riverine environment.
	\begin{align} 
		\label{reduced_system}
		\begin{cases}
			\vspace{0.2cm}
			&\displaystyle\frac{\partial u}{\partial t}= ({\mathsf{d}_1}{u_x}-{\alpha_1}{u})_x+ru{\left[1-\frac{u+v}{K}\right]},\ \displaystyle a<x<b,\ t>0, \\
			\vspace{0.2cm}
			&\displaystyle\frac{\partial v}{\partial t}= ({\mathsf{d}_2}{v_x}-{\alpha_2}{v})_x+rv{\left[1-\frac{u+v}{K}\right]},\ \displaystyle a<x<b,\ t>0,\\
			\vspace{0.2cm}
			&{\mathsf{d}_1}{u}_x(a,t)-{\alpha_1}u(a,t)=	{\mathsf{d}_1}{u}_x(b,t)-{\alpha_1}u(b,t)=0, \ \displaystyle t>0,\\
			\vspace{0.2cm}
			&{\mathsf{d}_2}{v}_x(a,t)-{\alpha_2}v(a,t)={\mathsf{d}_2}{v}_x(b,t)-{\alpha_2}v(b,t)=0,\ \displaystyle t>0,\\
			\vspace{0.2cm}
			&	u(x,0)={u_0}(x)\geq\not\equiv0,\ \displaystyle a<x<b,\\
			\vspace{0.2cm}
			&		v(x,0)=v_0(x)\geq\not\equiv0,\ \displaystyle a<x<b.
		\end{cases}
	\end{align}
The carrying capacity of the ecosystem is denoted by $K$. The first and second species densities at time $t$ and place $x$ are denoted by $u(x,t)$\ and $v(x,t)$, respectively. $t$, $(b-a)$ is the river's length, where $x=a$ denotes the upstream end (the water's surface) and $x=b$ denotes the downstream end (the bottom). The diffusion coefficients of the first and second species are $\mathsf{d}_1$ and $\mathsf{d}_2$, respectively, and the advection rates are $\alpha_1$ and $\alpha_2$. Given that $x=b$\ is designated as the downstream end, $\alpha$ ought to be positive. The existence of spatial heterogeneity in the environment is indicated by the intrinsic growth rate, or constant $r$\. Presume that $\mathsf{d}, r, a$, and $b$ are all positive constants.
%

In this study, we examined a reaction-diffusion-advection model designed to capture the dynamics of two competing species within a spatially heterogeneous riverine habitat. A key feature of the model is the inclusion of distinct diffusion and advection rates for each species, reflecting differences in their movement behaviors and interactions with the environment. Unlike the reduced model presented in \eqref{reduced_system}, which focuses exclusively on the effects of diffusion and advection on the population dynamics, the current study provides a more comprehensive analysis by incorporating the additional influence of harvesting. By integrating the effects of diffusion, advection, and harvesting simultaneously, we aim to explore the combined and individual roles of these factors in shaping the competitive dynamics, spatial distribution, and long-term persistence of the two species. This approach allows for a deeper understanding of how environmental flow, species mobility, and external exploitation interact to influence population stability and coexistence in ecological systems.
When harvesting is performed to a single or interacting species with various diffusive tactics, several scenarios might occur \cite{mk_ellena_harvesting,ellena_harvesting}. Previously, a single species with spatiotemporal harvesting was examined to investigate population persistence and extinction \cite{mk_ellena_harvesting}. In contrast, \cite{ellena_harvesting} focused on harvesting effort by utilizing a time-independent cropping function for two competitive populations. A regionally distributed harvesting map was used to investigate a non-homogeneous Gilpin-Ayala diffusive equation for single species \cite{bai_wang_harvesting}. In a current study, Adan et al. \cite{mk_adan} considered a competition model with the harvesting effect and examine that the most realistic harvesting scenario relates to the intrinsic growth rate. Additionally, Kamrujjaman and Zahan in \cite{mk_iz} demonstrate the impact of intrinsic growth rates and harvesting rates by examining a reaction-diffusion equation with resource-based diffusion strategies while taking into account two distinct scenarios.
	
Though a study on the harvesting effect for two species in \cite{mk_adan} reaction-diffusion model is efficient to provide an effective explanation of the actual ecological dynamics. But by taking into consideration the harvesting effect on the reaction-diffusion-advection (RDA) system on both species will give proper explanation.
	
We have evaluated the theoretical and numerical findings and looked at for very small harvesting rate how various advection to diffusion rates affect the system. Our findings reveal that under conditions of lower harvesting pressure, a species is more likely to outcompete its rival and establish dominance if it exhibits a slower advection rate coupled with a higher diffusion rate. This combination leads to a comparatively lower advection-to-diffusion ratio, providing the species with a strategic advantage in resource utilization and spatial distribution. Specifically, in a competitive environment, the ability of a species to diffuse more effectively allows it to explore and occupy a broader spatial area, thereby accessing more resources. Simultaneously, a reduced advection rate minimizes the risk of being swept away by environmental flows, enabling the species to maintain a stable presence in favorable regions. Thus, species characterized by a higher diffusion rate and a strictly lower advection rate have a significantly greater likelihood of persisting and thriving in the face of competition. These results highlight the critical role of movement dynamics in shaping the outcomes of interspecies competition and underline the importance of the advection-to-diffusion ratio as a key determinant of survival in ecological systems.
	
In our study, we consider a Lotla-Volterra two species reaction-diffusion-advection model with small harvesting effect.
The main findings of this paper are summarized as follows:
	\begin{enumerate}
		\item We establish the existence, uniqueness, and positivity of solutions for a reaction-diffusion-advection (RDA) model that incorporates the effects of harvesting within a heterogeneous advective environment. The model is analyzed under the condition that no-flux boundary constraints are imposed, which reflect the biological scenario where there is no exchange of individuals across the boundaries of the habitat. The existence of solutions ensures that the model accurately represents the biological dynamics over time, while uniqueness guarantees that the outcomes are well-defined and consistent for a given set of initial and boundary conditions. Positivity is a critical property in ecological models, as it ensures that population densities remain non-negative throughout the domain, reflecting realistic biological conditions. These mathematical properties form a foundational framework for analyzing the stability, dynamics, and long-term behavior of the population under the combined influences of spatial heterogeneity, advection, and harvesting. By addressing these properties rigorously, the study provides a robust theoretical basis for further exploration of ecological and management implications in advective environments.
		\item In the case of space-dependent carrying capacity with uneven diffusion and advection rate, two semi-trivial steady states are shown to be locally stable.
		\item Our analysis indicates that a species with a lower diffusion rate can consistently outcompete its rival if its advection rate, coupled with the harvesting rate, remains within the interval $(0,1]$. This outcome is particularly evident when the ratio of its advection rate to diffusion rate is smaller compared to the competing species. A lower advection-to-diffusion ratio allows the species to maintain a more stable presence within the environment, enabling it to access resources effectively while minimizing the dispersal disadvantage posed by excessive advection.
		
		Furthermore, the stability of the coexistence steady state is influenced by small perturbations in the advection rate of the second species. Minor changes to the advection dynamics of this competing species can tip the balance between competitive exclusion and coexistence, suggesting that the system's response to such perturbations is a key determinant of long-term population dynamics. This interplay highlights the delicate equilibrium within heterogeneous environments, where species movement behaviors, combined with external factors like harvesting, dictate outcomes ranging from dominance to coexistence. Understanding these nuanced interactions offers critical insights into ecological balance and resource management in advective systems.
		\item The numerical simulations presented in this study provide robust validation for the theoretical findings regarding species competition. These examples demonstrate a range of possible outcomes, including global asymptotic stability, coexistence, and extinction, under various parameter configurations and initial conditions. Global asymptotic stability is observed when one species consistently outcompetes the other, driving it to extinction while stabilizing its own population at an equilibrium state. Coexistence emerges as a viable outcome when both species are able to maintain positive population densities over time, indicating a balance in competitive pressures facilitated by differences in movement dynamics, harvesting rates, or environmental heterogeneity. Extinction scenarios occur when one or both species fail to sustain a viable population, often due to adverse combinations of advection, diffusion, or excessive harvesting.
		
		The numerical results not only confirm the analytical predictions but also highlight the sensitivity of the system to parameter variations. These simulations serve as a practical illustration of the complex interplay between species dynamics, spatial heterogeneity, and external interventions like harvesting, providing deeper insights into the conditions that promote stability, coexistence, or collapse within ecological systems.
	\end{enumerate}
The remainder of this paper is organized as follows:

In Section \ref{section_2}, we introduce and examine a two-species Lotka-Volterra reaction-diffusion-advection (RDA) model that incorporates harvesting effects. This section outlines the mathematical framework and governing equations that form the foundation of the study.
Section \ref{section_3} addresses the theoretical underpinnings of the model by establishing the existence, uniqueness, and positivity of its solutions. These properties ensure that the model accurately represents realistic ecological dynamics, with well-defined and biologically meaningful outcomes.
Preliminary results derived from the model are presented in Section \ref{section_4}, providing insights into the basic interactions and behaviors of the competing species under the modeled conditions.
In Section \ref{section_5}, we analyze the local stability of two semi-trivial steady states within a heterogeneous environment, highlighting the conditions under which each species can dominate in the absence of coexistence.
Section \ref{section_6} explores scenarios where the harvesting rate falls within the interval $(0,1]$. Here, we demonstrate the nonexistence of a coexistence steady state by employing specific and non-trivial assumptions regarding unequal diffusion and advection rates for the two species. This section builds on a monotone dynamical system framework and combines local stability analysis with the nonexistence of coexistence states to evolve the primary findings of this research.
Finally, Section \ref{section_7} provides numerical examples that validate and support the analytical results. These examples illustrate the practical implications of the model, showcasing the dynamics of species competition under varying conditions and confirming the theoretical conclusions.

This structured approach ensures a comprehensive exploration of the model and its implications, culminating in a detailed understanding of species interactions in advective environments with harvesting.

\section{Model Description}\label{section_2}
Suppose a Lotka-Volterra RDA model that includes harvesting for two competitive species with varying advection, diffusion, and harvesting rates, that is, $\mathsf{d}_1\neq \mathsf{d}_2, \  \alpha_1\neq\alpha_2\ \text{and}\ \mu_1 \neq \mu_2$.
Then the model under consideration becomes,
\begin{align}\label{system1}
	\begin{cases}
		\vspace{0.2cm}
		&\displaystyle\frac{\partial u}{\partial t}= \left({\mathsf{d}_1}{u_x}-{\alpha_1}{u}\right)_x+ru{\left[1-\frac{u+v}{K}\right]}- \mu_1 ru,\quad\displaystyle a<x<b,\ t>0 \\
		\vspace{0.2cm}
		&\displaystyle \frac{\partial v}{\partial t}= ({\mathsf{d}_2}{v_x}-{\alpha_2}{v})_x+rv{\left[1-\frac{u+v}{K}\right]}-\mu_2 rv,\quad\displaystyle a<x<b,\ t>0 \\
		\vspace{0.2cm}
		&{\mathsf{d}_1}{u}_x(a,t)-{\alpha_1}u(a,t)=	{\mathsf{d}_1}{u}_x(b,t)-{\alpha_1}u(b,t)=0, \quad\displaystyle t>0,\\
		\vspace{0.2cm}
		&{\mathsf{d}_2}{v}_x(a,t)-{\alpha_2}v(a,t)={\mathsf{d}_2}{v}_x(b,t)-{\alpha_2}v(b,t)=0,\quad\displaystyle t>0,\\
		\vspace{0.2cm}
		&	u(x,0)={u_0}(x)\geq,\not\equiv0,\ \displaystyle a<x<b\\
		\vspace{0.2cm}
		&		v(x,0)=v_0(x)\geq ,\not\equiv0, \quad \displaystyle a<x<b
	\end{cases}
\end{align}
The carrying capacity of the ecosystem is denoted by $K$. 
The diffusion coefficients of the first and second species are $\mathsf{d}_1$ and $\mathsf{d}_2$, respectively, and the advection rates are $\alpha_1$ and $\alpha_2$. Given that $x=b$\ is designated as the downstream end, $\alpha$ ought to be positive. The existence of spatial heterogeneity in the environment is indicated by the intrinsic growth rate, or constant $r$. Presume that $\mathsf{d}, r,\; a$, and $b$ are all positive constants. We will develop new and nontrivial techniques for equal harvesting rates (i.e., $\mu_1=\mu_2=\mu$) for the establishment of the non-existence result of any coexistence stable state, which is generally a challenging task in the study of monotone dynamical systems, to solve the developing mathematical challenges. But we consider the unequal harvesting effect on both species in our numerical analysis.
Now consider the model where both species are harvested and suppose the harvesting rates are proportional to the space dependent intrinsic growth rates. 

In river ecosystems, food sources (like algae or detritus) and nutrients are not evenly distributed. The reaction-diffusion-advection model \eqref{system1} can simulate how these resources move and concentrate in different areas of the river. Understanding how this impacts the feeding patterns of fish, invertebrates, and other species can help in designing conservation strategies that ensure a stable and accessible food supply for river animals.

Then the model \eqref{system1} under consideration $\mu_1=\mu_2=\mu$ reads as follows
\begin{align} 
	\label{harvest_system}
	\begin{cases}
		\vspace{0.2cm}
		&\displaystyle\frac{\partial u}{\partial t}= [\mathsf{d}_1  u_x-\alpha_1u]_x+ru{\left[1-\frac{u+v}{K}\right]}- \mu r u,\quad\displaystyle a<x<b,\ t>0 \\
		\vspace{0.2cm}
		&\displaystyle \frac{\partial v}{\partial t}= [\mathsf{d}_2 v_x-\alpha_2 v]_x+rv{\left[1-\frac{u+v}{K}\right]}-\mu r v,\quad\displaystyle a<x<b,\ t>0 \\
		\vspace{0.2cm}
		&\mathsf{d}_1u_x-\alpha_1u=\mathsf{d}_2v_x-\alpha_2v=0,\ \displaystyle x=a,b\\
		\vspace{.3cm}
		&	u(x,0)={u_0}(x)\geq,\not\equiv0,\ \displaystyle a<x<b\\
		\vspace{0.2cm}
		&		v(x,0)=v_0(x)\geq ,\not\equiv0, \quad \displaystyle a<x<b
	\end{cases}
\end{align}
Now modify the system in such a way that the system will be harvesting free. Rewrite the first equation of \eqref{harvest_system} as
\begin{align*}
	\vspace{0.2cm}
	\displaystyle\frac{\partial u}{\partial t}&= [\mathsf{d}_1 u_x-\alpha_1u]_x+ru{\left[1-\frac{u+v}{K}\right]}- \mu r u,\\
	\vspace{0.2cm}
	&= [\mathsf{d}_1 u_x-\alpha_1u]_x+ru{\left[1-\frac{u+v}{K}-\mu\right]}\\
	\vspace{0.2cm}
	&=[\mathsf{d}_1 u_x-\alpha_1u]_x+rr_1u{\left[1-\frac{u+v}{K_1}\right]}
\end{align*} 
If we assume $r_1=1-\mu\ \text{and}\ K_1=Kr_1$. Then we obtain
\begin{align*}
	\displaystyle\frac{\partial u}{\partial t}&=[\mathsf{d}_1 u_x-\alpha_1u]_x+rr_1u{\left[1-\frac{u+v}{K_1}\right]}.
\end{align*}
Similarly, the second equation can be written as 
\begin{align*}
	\displaystyle\frac{\partial v}{\partial t}&=[\mathsf{d}_2 v_x-\alpha_2v]_x+rr_1v {\left[1-\frac{u+v}{K_1}\right]},
\end{align*}
Finally, the system \eqref{harvest_system} can be written of the following equivalent form
\begin{align}
	\label{harvest_free}
	\begin{cases}
		\vspace{.3cm}
		&\displaystyle\frac{\partial u}{\partial t}=[\mathsf{d}_1 u_x-\alpha_1u]_x+rr_1u{\left[1-\frac{u+v}{K_1}\right]},\\
		\vspace{.3cm}
		&\displaystyle\frac{\partial v}{\partial t}=[\mathsf{d}_2 v_x-\alpha_2v]_x+rr_1v {\left[1-\frac{u+v}{K_1}\right]},\\
		\vspace{.3cm}
		&\mathsf{d}_1u_x-\alpha_1u=\mathsf{d}_2v_x-\alpha_2v=0,\ \displaystyle x=a,b\\
		\vspace{.3cm}
		&	u(x,0)={u_0}(x)\geq,\not\equiv0,\ \displaystyle a<x<b,\\
		&		v(x,0)=v_0(x)\geq ,\not\equiv0, \quad\displaystyle a<x<b.
	\end{cases}
\end{align}
The scenario where\ $\mathsf{d}_1> \mathsf{d}_2>0\ \text{and}\ \alpha_1>\alpha_2>0$ is the first one we address. To gain a deeper understanding of the dynamics at play, we examine the system under the following two specific conditions. These conditions are designed to provide clearer insights into the competitive interactions between species and the factors influencing their outcomes. Our observations reveal that a critical determinant of the species' behavior is the ratio of the advection rate to the diffusion rate. This ratio plays a pivotal role in shaping the spatial and temporal dynamics of the population.

A lower advection-to-diffusion ratio indicates that a species can effectively disperse over a larger spatial area while minimizing the effects of advection-driven displacement, thereby allowing it to stabilize its population in favorable regions. Conversely, a higher ratio implies that advection dominates, potentially leading to reduced control over spatial positioning and increased susceptibility to environmental flow. By analyzing these conditions, we aim to unravel how this ratio influences the persistence, dominance, or extinction of species within a heterogeneous advective environment. These insights provide a foundation for understanding the complex interplay of movement dynamics and external factors such as harvesting in shaping ecological outcomes.
\begin{gather}
	\label{c_1}
	\displaystyle\frac{\alpha_1}{\mathsf{d}_1}\geq 	\displaystyle\frac{\alpha_2}{\mathsf{d}_2} \quad \left(\text{equivalently\ }\displaystyle\frac{\alpha_1}{\alpha_2}\geq 	        \displaystyle \frac{\mathsf{d}_1}{\mathsf{d}_2} \right)
\end{gather}
and
\begin{gather}
	\label{c_2}
	\displaystyle\frac{\alpha_1}{\mathsf{d}_1}< 	\displaystyle\frac{\alpha_2}{\mathsf{d}_2} \quad \left(\text{equivalently\ }\displaystyle\frac{\alpha_1}{\alpha_2}<	\displaystyle\frac{\mathsf{d}_1}{\mathsf{d}_2} \right)
\end{gather}
\section{Existence, Uniqueness and Positivity of Solution}\label{section_3}
To determine the uniqueness result, take a look at the following model for a single species in the domain $\Omega$ with a homogeneous Neumann boundary condition.
\begin{align} 
	\label{harvest_system2}
	\begin{cases}
		\vspace{0.2cm}
		&\displaystyle\frac{\partial u}{\partial t}= [\mathsf{d}_1 u_x-\alpha_1 u]_x+rr_1 u{\left[1-\frac{u}{K_1}\right]},\quad\displaystyle x\in \Omega,\ t>0 \\
		\vspace{0.2cm}
		&\mathsf{d}_1u_x-\alpha_1u=0,\ \displaystyle x\in\partial \Omega\\
		\vspace{0.2cm}
		&	u(x,0)={u_0}(x)\geq,\not\equiv0,\ \displaystyle x\in \Omega
	\end{cases}
\end{align}
\begin{Lemma}
	\label{lm17}
	Let $u_0(x) \in \Omega,$ be the initial continuous function that is both non-negative and non-trivial, and $u_0(x)>0 $ in a non-empty open bounded sub-domain $ \Omega_1 \subset \Omega $. There is a unique, positive solution $ u(t, x) $ to the problem \eqref{harvest_system} for any $ t>0 $.
\end{Lemma}
\begin{proof}
	Let us assume
	$$\mathcal{G}_1(x,u)=\mathcal{G}_2(x,u)u=rr_1u\left(1-\frac{u}{K_1}\right)$$
	where
	$$\mathcal{G}_2(x,u)=rr_1\left(1-\frac{u}{K_1}\right)$$
	Then the system \eqref{harvest_system2} becomes
	\begin{align} 
		\label{harvest_system22}
		\begin{cases}
			\vspace{0.2cm}
			&\displaystyle\frac{\partial u}{\partial t}= [\mathsf{d}_1 u_x-\alpha_1 u]_x+\mathcal{G}_2u,\quad\displaystyle x\in \Omega,\ t>0 \\
			\vspace{0.2cm}
			&\mathsf{d}_1u_x-\alpha_1u=0,\ \displaystyle x\in\partial \Omega\\
			\vspace{0.2cm}
			&	u(x,0)={u_0}(x)\geq,\not\equiv0,\ \displaystyle x\in \Omega
		\end{cases}
	\end{align}
	In this instance, $\mathcal{G}_1(x,u)$ is Lipschitz in $u$ and is a measurable function in $x$ that is bounded if $u$ is limited to a bounded set, $\Omega$ is bounded, and $\partial\Omega$ belongs to class $C^{2+{\alpha}}(\Omega)$. The class $C^2 $ in $u $ contains the functions $\mathcal{G}_1(x,u)$ and $\mathcal{G}_2(x,u)$. For $ u > K_1 $, there exists $ K_1> 0 $ such that $\mathcal{G}_2(x,u)<0$. \eqref{harvest_system22}'s associated eigenvalue problem is shown as follows
	\begin{align}
		\label{ef}
		\begin{cases}
			&[\mathsf{d}_1\Phi_x-\alpha_1\Phi]_x+\mathcal{G}_2(x,0)=\lambda\Phi,\quad x\in \Omega,\\
			&\mathsf{d}_1\Phi_x-\alpha_1\Phi=0,\quad x\in\partial\Omega
		\end{cases}
	\end{align}
	We may write $ \mathcal{G}_1(x,u) = (\mathcal{G}_2(x,u) + \mathcal{G}_2(x,u)u) u $ based on the assumptions on $\mathcal{G}_1(x,u)$, if this problem's primary eigenvalue $\lambda$ is positive. For \eqref{ef}, let $\xi$ be an eigenfunction such that $\xi>0$ on $\Omega$. If $\epsilon>0$ is small enough,
	\begin{align*}
		\mathsf{d}_1\Delta(\epsilon\xi)-\alpha_1(\epsilon\xi)_x+\mathcal{G}_1(x,\epsilon\xi)&=\epsilon\left[\mathsf{d}_1\xi_x-\alpha_1\xi+\mathcal{G}_2(x,0)\xi\right]+\mathcal{H}(x,\epsilon\xi)\epsilon^2\xi^2\\
		&=\epsilon\lambda\xi+\mathcal{H}(x,\epsilon\xi)\epsilon^2\xi^2\\
		&=\epsilon\xi\{\lambda+\mathcal{H}(x,\epsilon\xi)\epsilon\xi\}
	\end{align*}
	Hence for $\epsilon > 0$ small, $\epsilon\xi$	 is a subsolution for the elliptic problem
	\begin{align}
		\label{ef2}
		\begin{cases}
			&[\mathsf{d}_1u_x-\alpha_1u]_x+\mathcal{G}_1(x,u)=\lambda u,\quad x\in \Omega,\\
			&\mathsf{d}_1u_x-\alpha_1u=0,\quad x\in\partial\Omega,\\
			&u(x,0)=u_0(x)\geq,\not\equiv0,\ \displaystyle x\in \Omega
		\end{cases}
	\end{align}
	corresponding to \eqref{harvest_system2}. At $t = 0 $, $\frac{\partial\underline{u}}{\partial t}>0$ on $\Omega$ if $\underline{u}(x,t)$ is a solution to \eqref{harvest_system2} with $\underline{u}(x,0) = \epsilon\xi$. Additionally, the general properties of sub-solutions and supersolutions suggest that $\underline{u}(x,t)$ is increasing in $t $. The minimal positive solution of \eqref{harvest_system2} is $u^\ast$, and since $ K_1 > u $ is a supersolution to \eqref{harvest_system2}, we must have $\underline{ u}(x,t)\uparrow u^\ast(x) $ as $ t\rightarrow\infty $. We can be certain that $u^\ast$ is minimal because $\xi$ will be a strict subsolution for all $ \epsilon> 0 $ sufficiently small. The strong maximum principle suggests that $ u(x,t) > 0 $ on $\overline{\Omega}$ for $ t > 0 $ if $ u(x,t) $ is a solution to \eqref{harvest_system2}, which is initially nonnegative and positive on an open subset of $\Omega$. This completes the evidence.
\end{proof}
\begin{Lemma}\label{lm18}
Equation \eqref{harvest_system2} has a unique equilibrium solution, $u_0^\ast(x)$ . The solution $u(t,x)$\ of \eqref{harvest_system2} then meets the condition $$\lim_{t\rightarrow\infty}u(t,x)=u^\ast(x)$$ uniformly for $x\in \overline{\Omega}$ for any initial solution $u_0(x)\ge 0,\ u_0(x)\not\equiv 0$.
\end{Lemma}
\begin{proof}
Assume that $\mathcal{G}_1(x,u) = \mathcal{G}_2(x,u)u $ with $\mathcal{G}(x,u)$ strictly decreasing in $u $ for $ u \ge 0 $, and that the hypotheses of Lemma \ref{lm17} are met. Then, for \eqref{harvest_system2}, the only positive equilibrium is the smallest positive equilibrium $u^\ast$. Since $u^\ast $ is minimum, we must have $u^{\ast\ast}> u^\ast $ somewhere on $\Omega$ if $u^{\ast\ast}$ is another positive equilibrium of \eqref{harvest_system2} with $u^{\ast\ast}\ne u^\ast$. $ u^\ast > 0 $ is a positive solution to \eqref{harvest_system2} since it is an equilibrium of
	\begin{align}
		\label{ef3}
		\begin{cases}
			&[\mathsf{d}_1\Phi_x-\alpha_1\Phi]_x+\mathcal{G}_2(x,u^\ast)=\lambda\Phi,\quad x\in \Omega,\\
			&\mathsf{d}_1\Phi_x-\alpha_1\Phi=0,\quad x\in\partial\Omega
		\end{cases}
	\end{align}
	for any $\Phi$, if we let $\lambda=0$, then $\lambda=0$ must be the principal eigenvalue of \eqref{ef3}. Similarly $ u^{\ast\ast}$ satisfies
	\begin{align}
		\label{ef4}
		\begin{cases}
			&[\mathsf{d}_1\Phi_x-\alpha_1\Phi]_x+\mathcal{G}_2(x,u^{\ast\ast})=\lambda\Phi,\quad x\in \Omega,\\
			&\mathsf{d}_1\Phi_x-\alpha_1\Phi=0,\quad x\in\partial\Omega
		\end{cases}
	\end{align}
	with $\lambda=0$ for any $\Phi$, so $\lambda=0$ in \eqref{ef4} also. But since $ \mathcal{G}_2(x,u) $ is strictly decreasing in $u $ and $ u^{\ast\ast}> u^\ast $ on at least a portion of $\Omega$, the principal eigenvalue in \eqref{ef4} needs to be smaller than the principal eigenvalue in \eqref{ef3}. The proof is therefore completed since we cannot have $\lambda=0$ in either \eqref{ef3} or \eqref{ef4}. Consequently, \eqref{harvest_system2} can only have the minimal equilibrium $u^\ast$.
\end{proof}
\begin{Lemma}
	\label{lm19}
	Let  $ v_0(x) > 0 $ in some non-empty open bounded sub-domain $ \Omega_1 \subset \Omega $, and $ v_0(x) \in \Omega,$ be the non-negative and non-trivial initial continuous function. The system \eqref{harvest_system} then has a unique, positive solution $ v(t, x) $ for any $ t > 0 $.
\end{Lemma}
\begin{Lemma}\label{lm20}
	The equation \eqref{harvest_system2}, has a unique equilibrium solution $v_0^\ast(x)$. The solution $v(t,x)$\ of \eqref{harvest_system2} then meets the condition $$\lim_{t\rightarrow\infty}v(t,x)=v^\ast(x)$$ uniformly for $x\in \overline{\Omega}$ for any initial solution $v_0(x)\ge 0,\ v_0(x)\not\equiv 0$.
\end{Lemma}
\begin{proof}
	The Lemmas \ref{lm19} and \ref{lm20} are analogous to Lemma \ref{lm17} and \ref{lm18}, so we can omit the proofs here. 
\end{proof}
Consider the stationary problem as
\begin{align}
	\label{harvest_stationary}
	\begin{cases}
		&[\mathsf{d}_1 u_x-\alpha_1u]_x+{r_1}ru{\left[1-\displaystyle\frac{u}{K_1}\right]}=0,\quad\displaystyle a<x<b,\ t>0 \\
		&\mathsf{d}_1u_x-\alpha_1u=0,\ \displaystyle x=a,b,\\
		&	u(x,0)={u_0}(x)\geq,\not\equiv0,\ \displaystyle a<x<b
	\end{cases}
\end{align}
\begin{Lemma}\label{L22}
	Suppose that $K_1(x)\not\equiv \text{const}$ and $u^\ast(x)$ is the positive solution of \eqref{harvest_stationary}, then
	$$\displaystyle\int_a^b{r(x)K_1\left(1-\frac{u^\ast(x)}{K_1}\right)}dx>0.$$
\end{Lemma}
\begin{proof}
	Since $K(x)$\ and\ $u^\ast(x)$\ are strictly positive, then integrating the first equation of \eqref{harvest_stationary} over the domain, we get
	\begin{align*}
		&\displaystyle\int_a^b{[\mathsf{d}_1 u_x^{\ast}(x)-\alpha_1u^\ast(x)]_x+{r_1}r(x)u^\ast(x){\left[1-\frac{u^\ast(x)}{K_1}\right]}}dx=0
	\end{align*}
	By applying homogeneous Neumann boundary condition we obtain
	\begin{align*}
		&\displaystyle\int_a^b{{r_1}r(x)u^\ast(x){\left[1-\frac{u^\ast(x)}{K_1(x)}\right]}}dx=0\\
		&\Rightarrow\displaystyle\int_a^b{r(x)(u^\ast(x)-K_1+K_1){\left[1-\frac{u^\ast(x)}{K_1}\right]}}dx=0\\
		&\Rightarrow\displaystyle\int_a^b{r(x)(u^\ast(x)-K_1){\left[1-\frac{u^\ast(x)}{K_1}\right]}}dx+\displaystyle\int_a^b{r(x)K_1{\left[1-\frac{u^\ast(x)}{K_1}\right]}}dx=0\\
		&\Rightarrow\displaystyle\int_a^b{r(x)K_1{\left[1-\frac{u^\ast(x)}{K_1}\right]^2}}dx=\displaystyle\int_a^b{r(x)K_1{\left[1-\frac{u^\ast(x)}{K_1}\right]}}dx
	\end{align*}
	Clearly for $u^\ast(x)\ne K_1,\ \displaystyle\int_a^b{r(x)K_1{\left[1-\frac{u^\ast(x)}{K_1}\right]^2}}dx>0$. Since $u^\ast(x)=K_1\not\equiv\text{const}$, then $u^\ast(x)$ is not a solution of \eqref{harvest_stationary}. Therefore,
	$$\displaystyle\int_a^b{r(x)K_1{\left[1-\frac{u^\ast(x)}{K_1}\right]}}dx>0$$
\end{proof}
Consider the stationary problem as
\begin{align}
	\label{harvest_stationary2}
	\begin{cases}
		&[\mathsf{d}_2 v_x-\alpha_2 v]_x+{r_1}rv{\left[1-\frac{v}{K_1}\right]}=0,\quad\displaystyle a<x<b,\ t>0 \\
		&\mathsf{d}_2v_x-\alpha_2v=0,\ \displaystyle x=a,b\\
		&	v(x,0)={v_0}(x)\geq,\not\equiv0,\ \displaystyle x\in \Omega
	\end{cases}
\end{align}
\begin{Lemma}
	Suppose that $K_1(x)\not\equiv \text{const}$ and $v^\ast(x)$ is the positive solution of \eqref{harvest_stationary2}, then
	$$\displaystyle\int_a^b{r(x)K_1\left(1-\frac{v^\ast(x)}{K_1}\right)}dx>0.$$
\end{Lemma}
\begin{proof}
	The proof is analogous to the previous lemma \ref{L22} so we can omit it here.
\end{proof}
\begin{Th}
	For $\mu\in [0,1)$ and any $u_0(x),v_0(x)\in C(a,b)$, the system \eqref{harvest_system} and \eqref{harvest_free} have a unique solution $ (u,v)$. Furthermore, for every $ t > 0 $, $ u(t,x) > 0 $ and $ v(t,x) > 0 $ if both starting conditions $u_0$ and $v_0$ are non-negative and non-trivial.
\end{Th}
\begin{proof}
	According to \cite{e3}, the system \eqref{harvest_system} has a non trivial time dependent solution. Consider 
	$$\gamma_u=\max\{\sup_{a<x<b}u_0(x),K(x)\},~~\gamma_v=\max\{\sup_{a<x<b}v_0(x),K(x)\}$$
	Suppose that 
	\begin{align}
		&f_1(t,x,u,v)=rr_1u{\left[1-\frac{u+v}{K_1}\right]}\\
		& f_2(t,x,u,v)=rr_1v{\left[1-\frac{u+v}{K_1}\right]}
	\end{align}
	and define
	$$P^\ast\equiv \{(\check{u},\check{v})\in C([0,\infty)\times [a,b]):0\le \check{u}\le \gamma_u,0\le \check{v}\le \gamma_v\}$$
	The following requirements are met since the functions $f_1$\ and\ $f_2$\ are quasi-monotone non-increasing Lipshitz functions in $P^\ast$ \cite{book}. 
	\begin{align}
		&f_1(t,x,\gamma_u,0)\le 0\le f_1(t,x,0,\gamma_v)\\
		&f_2(t,x,0,\gamma_v)\le 0\le f_2(t,x,\gamma_u,0)
	\end{align}
	We then consider the class of continuous functions $C([0,\infty)\times [a,b])$ on $[0,\infty)\times [a,b]$\ such that $$(u_0,v_0)\in P^\ast $$ and the initial solution $u_0(x),v_0(x)$. All $(t,x)\in [0,\infty)\times [a,b]$ then have a unique solution $(u(t,x),v(t,x))$\ of \eqref{harvest_system} in $P^\ast$. As a result, the solution $(u(t, x), v(t, x))$ is unique and positive.
\end{proof}
\section{Auxiliary Results}\label{section_4}
In this section, we present a set of fundamental lemmas that will serve as essential tools for the subsequent analysis. These lemmas are carefully derived to address key aspects of the system's behavior and provide the mathematical foundation for our results in later sections.

As established in the preceding section, the population dynamics of the system \eqref{harvest_free} are intricately tied to its steady-state behavior. The steady states not only dictate the long-term outcomes of the population interactions but also play a crucial role in determining stability and persistence within the system. To gain deeper insights into these dynamics, we now turn our attention to the stationary problem associated with the system. This involves analyzing the system under equilibrium conditions, where time-dependent changes vanish, allowing us to focus solely on the spatial and parametric factors influencing the steady states.

By addressing the stationary problem, we aim to elucidate how diffusion, advection, and harvesting interplay to shape the system's equilibrium and provide a framework for understanding the conditions that lead to coexistence, dominance, or extinction of species. These foundational results will be instrumental in advancing our understanding of the system's behavior and its response to various ecological and environmental factors.
\begin{align}
	\label{stationary}
	\begin{cases}
		\vspace{0.2cm}
		&[\mathsf{d} \Psi_x-\alpha \Psi]_x+rr_1\Psi\left[1-\frac{\Psi}{K_1}\right]=0, \quad a<x<b,\\
		\vspace{0.2cm}
		& \mathsf{d} \Psi_x(a)-\alpha \Psi(a)=0\\
		\vspace{0.2cm}
		&\mathsf{d} \Psi_x(b)-\alpha \Psi(b)=0\\
	\end{cases}
\end{align} 
In this instance, $\mathsf{d}>0,\ r_1\in[0,1)\ \text{and}\ \alpha\in \mathbb{R}$.
\subsection{Existence of Semi-Trivial Steady State}
The following lemma establishes and guarantees the existence of two distinct semi-trivial steady states, denoted as $(\check{u},0)$ and 
$(0,\check{v})$, for the stationary problem \eqref{stationary}. These steady states represent scenarios where one species dominates while the other is absent, highlighting critical aspects of the competitive dynamics within the system.

Specifically, the steady state $(\check{u},0)$ corresponds to the scenario where the first species achieves equilibrium in the absence of the second species, while 
$(0,\check{v})$ represents the reverse situation, with the second species stabilizing without the presence of the first. The existence of these semi-trivial steady states is foundational, as it reflects the capacity of each species to survive independently under the given environmental conditions and parameter configurations.

This lemma not only confirms the mathematical feasibility of these steady states but also provides insights into the ecological implications, such as the conditions under which one species might dominate and exclude the other. These results serve as a precursor to further investigations into the stability and feasibility of coexistence or competitive exclusion within the system. By rigorously establishing the existence of these steady states, we set the stage for a comprehensive exploration of the interplay between diffusion, advection, and harvesting in shaping population outcomes.

\begin{Lemma}
	\label{lma1}
	There is a unique positive solution $\displaystyle\check\zeta,$ for every $K_1=\text{const}, \mu\in [0,1), \mathsf d>0$ and $\alpha\in \mathbb{R}$ in the system \eqref{stationary} that satisfies
	\begin{enumerate}
		\item[(a)] $\check{\beta}_x\equiv 0$ in $[a,b]$\ if $\alpha=0;$
		\item[(b)] $0< \displaystyle\frac{\check{\beta}_x}{\beta}<\displaystyle\frac{\alpha}{\mathsf{d}}$ in $[a,b]$\ if $\alpha>0;$
		\item[(c)] $0> \displaystyle\frac{\check{\beta}_x}{\beta}>\displaystyle\frac{\alpha}{\mathsf{d}}$ in $[a,b]$\ if $\alpha<0;$
	\end{enumerate}
\end{Lemma}
\begin{proof}
	When $\alpha=0,\ \check{\beta}=r$ in $(a,b)$, and therefore $\check{\beta}_x\equiv 0$, the proof of assertion (a) is straightforward.
	Let's examine the transformation to demonstrate statement (b)
	$$\mathcal{T}:=\displaystyle\frac{\check{\beta}_x}{\check{\beta}}$$
	We can derive \cite{main}
	\begin{align}
		\begin{cases}
			&-\mathsf{d}\Delta \mathcal{T}+[\alpha-2\mathsf{d}\mathcal{T}]\mathcal{T}_x+u\mathcal{T}=0, \quad a<x<b,\\
			&\mathcal{T}(a)= \mathcal{T}(b)=\frac{\alpha}{\mathsf{d}},\\
		\end{cases}
	\end{align}
	Then, applying the maximum principle at our disposal,
	$$0<\mathcal{T}<\displaystyle \frac{\alpha}{\mathsf{d}}\Rightarrow 0<\displaystyle\frac{\check{\beta}_x}{\check{\beta}}<\displaystyle \frac{\alpha}{\mathsf{d}} $$
	This results in statement (b). Likewise, we can derive assertion (c).
\end{proof}
\section{Local Stability of Semi-Trivial Steady States}
\label{section_5}
The eigenvalue problem will now be examined in order to examine the local stability of the semi-trivial steady state.
\begin{align}
	\label{eigen_eq2}
	\begin{cases}
		&[\mathsf{d}{\Phi_1}_x-\alpha\Phi_1]_x+[p(x)+\lambda_1]\Phi_1=0, \quad a<x<b\\
		& {\Phi_1}_x(a,t)-\alpha\Phi_1(a,t)={\Phi_1}_x(b,t)-\alpha\Phi_1(b,t)=0,
	\end{cases}
\end{align}
We shall assume that $p(x)=rr_1(1-\frac{\check{u}}{K_1})\in {C^1}([a,b])$ in this point. The first eigenvalue-eigenfunction pair of the equation \eqref{eigen_eq2} is represented as $(\lambda_1,\Phi_1)$. Keep in mind that $\lambda_1$ is simple and $\Phi_1$ is strictly positive in the domain.\\
\begin{Lemma}
	\label{lma2}
	The statements listed below are true 
	\begin{enumerate} 
		\item [(a)] In $(a,b)$, $\mathsf{d}{\Phi_1}_x<\alpha\Phi_1$, provided $p_x\leq,\not\equiv0$\ in $[a,b] $;
		 \item[(b)]  In $(a,b)$, $\mathsf{d}{\Phi_1}_x>\alpha\Phi_1$, provided $p_x\geq,\not\equiv0$\ in $[a,b] $;
		  \end{enumerate}
\end{Lemma}
\begin{proof}
	Suppose the transformation
	$$\mathcal{T}:=\displaystyle\frac{{\Phi_1}_x}{\Phi_1}$$
	We can derive \cite{main}
	\begin{align}
		\begin{cases}
			&-\mathsf{d}\Delta \mathcal{T}+[\alpha-2\mathsf{d}\mathcal{T}]\mathcal{T}_x+u\mathcal{T}=0, \quad x\in(a,b),\\
			&\mathcal{T}(a)= \mathcal{T}(b)=\frac{\alpha}{\mathsf{d}},\\
		\end{cases}
	\end{align}
	The aforementioned lemma can then be proved using the maximal principle.
\end{proof}

\begin{Lemma} \cite{main}
	\label{lma4}
	Suppose that $\mathsf{d}_1>\mathsf{d}_2>0$\ and $\alpha_1>\alpha_2>0$. Then
	$$\displaystyle\frac{\alpha_1}{\mathsf{d}_1}\geq\displaystyle\frac{\alpha_2}{\mathsf{d}_2}\Leftrightarrow\displaystyle\frac{\alpha_1-\alpha_2}{\mathsf{d}_1-\mathsf{d}_2}\geq\displaystyle\frac{\alpha_1}{\mathsf{d}_1}\Leftrightarrow\displaystyle\frac{\alpha_1-\alpha_2}{\mathsf{d}_1-\mathsf{d}_2}\geq\displaystyle\frac{\alpha_2}{\mathsf{d}_2}$$ 
\end{Lemma}

\begin{Lemma}
	\label{Lm_23}
	Assume that $K\not\equiv\text{const},\ \mu\in [0,1),\ \mathsf{d}_1,\mathsf{d}_2>0,\ \text{and}\ \alpha_1,\alpha_2\in \mathbb{R}$. Then
	\begin{align}
		\label{h1}
		&{\eta_1^2}-{\eta_1^1}=\displaystyle \frac{\displaystyle \int_{a}^{b}{[(\mathsf{d}_2-\mathsf{d}_1){\zeta_1^1}_{x}+(\alpha_1-\alpha_2)\zeta_1^1]\cdot[e^{-\frac{\alpha_2}{\mathsf{d}_2}x}\zeta_1^2]_x} dx}{\displaystyle \int_{a}^{b}{e^{-\frac{\alpha_2}{\mathsf{d}_2}x}\cdot \zeta_1^1 \cdot \zeta_1^2}\ dx}		
	\end{align}
	where $\eta_1$ and $\zeta_1$ are dependent on both $\mathsf{d}$ and $\alpha$. Also $(\eta_1^1,\zeta_1^1)\equiv \left(\eta_1(\mathsf{d}_1,\alpha_1),\zeta_1(\mathsf{d}_1,\alpha_1)\right)$ and $(\eta_1^2,\zeta_1^2)\equiv\left(\eta_1(\mathsf{d}_2,\alpha_2),\zeta_1(\mathsf{d}_2,\alpha_2)\right).$
\end{Lemma}
\begin{proof}
	The linear eigenvalue problem of \eqref{harvest_free} will be used to investigate the local stability of a semi-trivial steady state.
	\begin{align}
		\begin{cases}
			\vspace{.2cm}
			&	[\mathsf{d}_1\zeta_x-\alpha_1\zeta]_x+p(x)\zeta+\eta\zeta=0,\quad a<x<b\\
			\vspace{.2cm}
			&	\mathsf{d}_1\zeta_x(a,t)-\alpha_1\zeta(a,t)=\mathsf{d}_1\zeta_x(b,t)-\alpha_1\zeta(b,t)=0,\ t>0\\
		\end{cases}
	\end{align}
	which can be rewrite as
	\begin{align}
		\begin{cases}
			\vspace{.2cm}
			&	[\mathsf{d}_2\zeta_{x}-\alpha_2\zeta]_x+p(x)\zeta+\eta\zeta=\left[(\mathsf{d}_2-\mathsf{d}_1)\zeta_{x}+(\alpha_1-\alpha_2)\zeta\right]_x,\quad a<x<b\\
			\vspace{.2cm}
			&	\mathsf{d}_2\zeta_x(a,t)-\alpha_2\zeta(a,t)=\mathsf{d}_2\zeta_{x}(b,t)-\alpha_2\zeta(b,t)=0,\qquad t>0\\
		\end{cases}
	\end{align}
	where we assume  $p(x)=rr_1(1-\frac{\check{u}}{K_1})\in {C^1}([a,b]).$
	Now rearrange the equations of $(\eta_1^1,\zeta_1^1)$\ and $(\eta_1^2,\zeta_1^2)$\ as follows:
	\begin{align*}
		\begin{cases}
			\vspace{.2cm}
			&	[\mathsf{d}_2{\zeta_1^1}_x-\alpha_2{\zeta_1^1}]_x+[p(x)+\eta_1^1]\zeta_1^1=[(\mathsf{d}_2-\mathsf{d}_1){\zeta_1^1}_x+(\alpha_1-\alpha_2){\zeta_1^1}]_x,\quad a<x<b\\	\vspace{.2cm}
			&	[\mathsf{d}_2{\zeta_1^2}_x-\alpha_2{\zeta_1^2}]_x+[p(x)+\eta_1^2]\zeta_1^2=0\\
			\vspace{.2cm}
			&	\mathsf{d}_2{\zeta_1^1}_x(a,t)-\alpha_2{\zeta_1^1}(a,t)=[\mathsf{d}_2-\mathsf{d}_1]{\zeta_1^1}_x(a,t)+[\alpha_1-\alpha_2]\zeta_1^1(a,t),\quad t>0\\
			\vspace{.2cm}
			&	\mathsf{d}_2{\zeta_1^1}_x(b,t)-\alpha_2\zeta_1^1(b,t)=[\mathsf{d}_2-\mathsf{d}_1]{\zeta_1^1}_x(b,t)+[\alpha_1-\alpha_2]\zeta_1^1(b,t),\quad t>0\\
			\vspace{.2cm}
			&	\mathsf{d}_2{\zeta_1^2}_x(a,t)-\alpha_2\zeta_1^2(a,t)=\mathsf{d}_2{\zeta_1^2}_x(b,t)-\alpha_2\zeta_1^2(b,t)
		\end{cases}
	\end{align*}
	Multiply the first equation by $e^{-\frac{\alpha_2}{\mathsf{d}_2}x}\zeta_1^2$\ and then integrate over $[a,b],$\ we obtain
	\vspace{-.2cm}
	\begin{multline*}
		\displaystyle \int_{a}^{b}{[\mathsf{d}_2{\zeta_1^1}_x-\alpha_2{\zeta_1^1}]_xe^{-\frac{\alpha_2}{\mathsf{d}_2}x}\zeta_1^2+
			[p(x)+\eta_1^1]\zeta_1^1e^{-\frac{\alpha_2}{\mathsf{d}_2}x}\zeta_1^2} dx=\\ 
		\displaystyle \int_{a}^{b}{\left[(\mathsf{d}_2-\mathsf{d}_1){\zeta_1^1}_x+(\alpha_1-\alpha_2){\zeta_1^1}\right]_xe^{-\frac{\alpha_2}{\mathsf{d}_2}x}\zeta_1^2	} dx
	\end{multline*}
	\vspace{-.4cm}
	\begin{multline}
		\label{int_11}
		\Rightarrow	\frac{1}{\mathsf{d}_2}\displaystyle \int_{a}^{b}{[\mathsf{d}_2{\zeta_1^1}_x-\alpha_2{\zeta_1^1}]}\cdot e^{-\frac{\alpha_2}{\mathsf{d}_2}x}\cdot[\mathsf{d}_2{\zeta_1^2}_x-\alpha_2{\zeta_1^2}] dx=
		\displaystyle\int_{a}^{b}{[p(x)+\eta_1^1]\zeta_1^1e^{-\frac{\alpha_2}{\mathsf{d}_2}x}\zeta_1^2} dx+\\ \displaystyle\int_{a}^{b}{\left[(\mathsf{d}_2-\mathsf{d}_1){\zeta_1^1}_x+(\alpha_1-\alpha_2){\zeta_1^1}\right]}\cdot\left[e^{-\frac{\alpha_2}{\mathsf{d}_2}x}\zeta_1^2\right]_x dx	
	\end{multline}
	In the same way, multiply the second equation by $e^{-\frac{\alpha_2}{\mathsf{d}_2}x}\zeta_1^1$\ and then integrate over $[a,b]$, we get
	\begin{equation*}
		\displaystyle \int_{a}^{b}{	[\mathsf{d}_2{\zeta_1^2}_x-\alpha_2{\zeta_1^2}]_xe^{-\frac{\alpha_2}{\mathsf{d}_2}x}\zeta_1^1+
			[p(x)+\eta_1^2]\zeta_1^2e^{-\frac{\alpha_2}{\mathsf{d}_2}x}\zeta_1^1} dx=0
	\end{equation*}
	\vspace{-.6cm}
	\begin{multline}
		\label{int_22}
		\Rightarrow	\displaystyle\frac{1}{\mathsf{d}_2}\displaystyle \int_{a}^{b}{[\mathsf{d}_2{\zeta_1^1}_x-\alpha_2{\zeta_1^1}]}\cdot e^{-\frac{\alpha_2}{\mathsf{d}_2}x}\cdot[\mathsf{d}_2{\zeta_1^2}_x-\alpha_2{\zeta_1^2}] dx=
		\displaystyle\int_{a}^{b}{[p(x)+\eta_1^2]\zeta_1^2e^{-\frac{\alpha_2}{\mathsf{d}_2}x}\zeta_1^1} dx	
	\end{multline}
	From the equations \eqref{int_11} and \eqref{int_22} it follows that
	\begin{align*}
		&[{\eta_1^2}-{\eta_1^1}]\cdot \displaystyle \int_{a}^{b}{e^{-\frac{\alpha_2}{\mathsf{d}_2}x}\cdot \zeta_1^1 \cdot \zeta_1^2}\ dx=\displaystyle \int_{a}^{b}{[(\mathsf{d}_2-\mathsf{d}_1){\zeta_1^1}_x+(\alpha_1-\alpha_2)\zeta_1^1]\cdot\left[e^{-\frac{\alpha_2}{\mathsf{d}_2}x}\zeta_1^2\right]_x} dx\nonumber\\
		&\Rightarrow {\eta_1^2}-{\eta_1^1}=\displaystyle \frac{\displaystyle \int_{a}^{b}{\left[(\mathsf{d}_2-\mathsf{d}_1){\zeta_1^1}_x+(\alpha_1-\alpha_2)\zeta_1^1\right]\cdot\left[e^{-\frac{\alpha_2}{\mathsf{d}_2}x}\zeta_1^2\right]_x} dx}{\displaystyle \int_{a}^{b}{e^{-\frac{\alpha_2}{\mathsf{d}_2}x}\cdot \zeta_1^1 \cdot \zeta_1^2}\ dx}		
	\end{align*}
\end{proof}

\begin{Lemma}
	\label{lm24}
	It is assumed that $K_1\not\equiv\text{const},\ \mu\in[0,1),\ \mathsf{d}_1>\mathsf{d}_2>0,\ \alpha_1>\alpha_2>0$\ and \eqref{c_1} are true. This indicates that $(\hat{u},0)$ is unstable.
\end{Lemma}
\begin{proof}
	To prove the lemma, recall the problem 
	\begin{align*}
		\begin{cases}
			\vspace{.2cm}
			&	({\mathsf{d}_1}{ \hat{u}_x}-{\alpha_1} \hat{u})_x+rr_1\hat{u}{[1-\frac{\hat{u}}{K_1}]}=0,\quad a<x<b\\
			\vspace{0.2cm}
			&{\mathsf{d}_1}{ \hat{u}}_x(a,t)-{\alpha_1} \hat{u}(a,t)=	{\mathsf{d}_1}{ \hat{u}}_x(b,t)-{\alpha_1} \hat{u}(b,t)=0, \quad\displaystyle t>0,\\
		\end{cases}
	\end{align*}
	and we write the eigenvalue problem as
	\begin{align*}
		\begin{cases}
			\vspace{.2cm}
			&	({\mathsf{d}_2}{\Phi_1}_x-{\alpha_2}{\Phi_1})_x+rr_1\Phi_1{[1-\frac{\hat{u}}{K_1}]}+\kappa_1
			\Phi_1=0,\quad a<x<b\\
			\vspace{0.2cm}
			&{\mathsf{d}_2}{\Phi_1}_x(a,t)-{\alpha_1}\Phi_1(a,t)=	{\mathsf{d}_2}{\Phi_1}_x(b,t)-{\alpha_1}\Phi_1(b,t)=0, \quad\displaystyle t>0,\\	
		\end{cases}
	\end{align*}
	Then according to the expression \eqref{h1} in Lemma \ref{Lm_23}, we have
	\begin{align*}
		\vspace{0.2cm}
		\kappa_1 &=\displaystyle\frac{\displaystyle\int_{a}^{b}{[(\mathsf{d}_2-\mathsf{d}_1)\hat{u}_x+(\alpha_1-\alpha_2)\hat{u}]\cdot\left[e^{-\frac{\alpha_2}{\mathsf{d}_2}x}\Phi_1\right]_x}dx}{\displaystyle\int_{a}^{b}{e^{-\frac{\alpha_2}{\mathsf{d}_2}x}\cdot\hat{u}\cdot\Phi_1 dx}}\\
		\vspace{0.2cm}
		&=\displaystyle\frac{\displaystyle\int_{a}^{b}{[(\mathsf{d}_2-\mathsf{d}_1)\hat{u}_x+(\alpha_1-\alpha_2)\hat{u}]\cdot \left[e^{-\frac{\alpha_2}{\mathsf{d}_2}x}{\Phi_1}_x-e^{-\frac{\alpha_2}{\mathsf{d}_2}x}\Phi_1\cdot\displaystyle\frac{\alpha_2}{\mathsf{d}_2}\right]
			}dx}{\displaystyle\int_{a}^{b}{e^{-\frac{\alpha_2}{\mathsf{d}_2}x}\cdot\hat{u}\cdot\Phi_1 dx}}\\
		\vspace{0.2cm}
		&=-\displaystyle\frac{\displaystyle\int_{a}^{b}{(\mathsf{d}_1-\mathsf{d}_2)\cdot\left[\displaystyle\frac{\hat{u}_x}{\hat{u}}-\displaystyle\frac{\alpha_1-\alpha_2}{\mathsf{d}_1-\mathsf{d}_2}	\right]\cdot
				\left[\displaystyle\frac{{\Phi_1}_x}{\Phi_1}-\frac{\alpha_2}{\mathsf{d}_2}	\right]\cdot e^{-\frac{\alpha_2}{\mathsf{d}_2}x}\cdot \hat{u}\cdot \Phi_1}dx}{\displaystyle\int_{a}^{b}{e^{-\frac{\alpha_2}{\mathsf{d}_2}x}\cdot\hat{u}\cdot\Phi_1 dx}}
	\end{align*}
	Given that $\alpha_1>0$, we can obtain from assertion (b) in Lemma \ref{lma1} that
	$$0<\displaystyle\frac{\hat{u}_x}{\hat{u}}<\displaystyle\frac{\alpha_1}{\mathsf{d}_1}$$
	From Lemma \ref{lma4}, one can write
	$$0<\displaystyle\frac{\hat{u}_x}{\hat{u}}<\displaystyle\frac{\alpha_1}{\mathsf{d}_1}\leq \displaystyle\frac{\alpha_1-\alpha_2}{\mathsf{d}_1-\mathsf{d}_2}$$
	This implies
	$$\displaystyle\frac{\hat{u}_x}{\hat{u}}-\displaystyle\frac{\alpha_1-\alpha_2}{\mathsf{d}_1-\mathsf{d}_2}<0$$
	Also, since $p(x)=rr_1(1-\frac{\hat{u}}{K_1})$, so
	$$p_x=-\frac{rr_1}{K_1}\hat{u}_x$$
	 In $[a,b]$, it is strictly negative. According to Lemma \ref{lma2}'s statement (b),
	$$\displaystyle\frac{{\Phi_1}_x}{\Phi_1}-\displaystyle\frac{\alpha_2}{\mathsf{d}_2}<0$$
   It is evident from combining all of the inequalities that $\kappa_1<0$, resulting in the instability of $(\hat{u},0)$. 
\end{proof}
\begin{Lemma}
	\label{lm25}
	It is assumed that $K_1\not\equiv\text{const},\ \mu\in[0,1),\ \mathsf{d}_1>\mathsf{d}_2>0,\ \alpha_1>\alpha_2>0$\ and \eqref{c_1} are true. This indicates that $(0,\hat{v})$\ is stable.
\end{Lemma}
\begin{proof}
	To prove this,  we write the eigenvalue problem as
	\begin{align*}
		\begin{cases}
			\vspace{.2cm}
			&	({\mathsf{d}_2}{\beta_1}_x-{\alpha_2}{\beta_1})_x+rr_1\hat{v}{[1-\frac{\hat{v}}{K_1}]}+\tau_1
			\beta_1=0,\quad a<x<b\\
			\vspace{0.2cm}
			&{\mathsf{d}_2}{\beta_1}_x(a,t)-{\alpha_1}\beta_1(a,t)=	{\mathsf{d}_2}{\beta_1}_x(b,t)-{\alpha_1}\beta_1(b,t)=0, \quad\displaystyle t>0,\\	
		\end{cases}
	\end{align*}
The expression \eqref{h1} of Lemma \ref{Lm_23} then gives us
	\begin{align*}
		\vspace{0.2cm}
		\tau_1&=\displaystyle\frac{\displaystyle\int_{a}^{b}{-[(\mathsf{d}_2-\mathsf{d}_1) \hat{v}_x+(\alpha_1-\alpha_2)\hat{v}]\cdot\left[e^{-\frac{\alpha_1}{\mathsf{d}_1}x}\beta_1\right]_x}dx}{\displaystyle\int_{a}^{b}{e^{-\frac{\alpha_1}{\mathsf{d}_1}x}\cdot\hat{v}\cdot\beta_1 dx}}\\
		\vspace{0.2cm}
		&=\displaystyle\frac{\displaystyle\int_{a}^{b}{-[(\mathsf{d}_2-\mathsf{d}_1)\hat{v}_x+(\alpha_1-\alpha_2)\hat{v}]\cdot \left[e^{-\frac{\alpha_1}{\mathsf{d}_1}x} {\beta_1}_x-e^{-\frac{\alpha_1}{\mathsf{d}_1}x}\beta_1\cdot\frac{\alpha_1}{\mathsf{d}_1}\right]
			}dx}{\displaystyle\int_{a}^{b}{e^{-\frac{\alpha_1}{\mathsf{d}_1}x}\cdot\hat{v}\cdot\beta_1 dx}}\\
		\vspace{0.2cm}
		&=\displaystyle\frac{\displaystyle\int_{a}^{b}{(\mathsf{d}_1-\mathsf{d}_2)\cdot\left[\displaystyle\frac{\hat{v}_x}{\hat{v}}-\displaystyle\frac{\alpha_1-\alpha_2}{\mathsf{d}_1-\mathsf{d}_2}	\right]\cdot
				\left[\displaystyle\frac{{\beta_1}_x}{\beta_1}-\displaystyle\frac{\alpha_1}{\mathsf{d}_1}	\right]\cdot e^{-\frac{\alpha_1}{\mathsf{d}_1}x}\cdot \hat{v}\cdot \beta_1}dx}{\displaystyle\int_{a}^{b}{e^{-\frac{\alpha_1}{\mathsf{d}_1}x}\cdot\hat{v}\cdot\beta_1 dx}}
	\end{align*}
	Given that $\alpha_1>0$, we can obtain from assertion (b) in Lemma \ref{lma1} that
	$$0<\frac{\hat{v}_x}{\hat{v}}<\displaystyle\frac{\alpha_1}{\mathsf{d}_1}$$
	From Lemma \ref{lma4}, one can write
	$$0<\displaystyle\frac{\hat{v}_x}{\hat{v}}<\displaystyle\frac{\alpha_1}{\mathsf{d}_1}\leq \displaystyle\frac{\alpha_1-\alpha_2}{\mathsf{d}_1-\mathsf{d}_2}$$
	This implies
	$$\displaystyle\frac{\hat{v}_x}{\hat{v}}-\displaystyle\frac{\alpha_1-\alpha_2}{\mathsf{d}_1-\mathsf{d}_2}<0$$
	Also, since $p(x)=rr_1(1-\frac{\hat{u}}{K_1})$, so
	$$p_x=-\frac{rr_1}{K_1}\hat{u}_x$$
	In $[a,b]$, it is strictly negative. According to Lemma \ref{lma2}'s statement (b),
	$$\displaystyle\frac{{\beta_1}_x}{\beta_1}-\frac{\alpha_1}{\mathsf{d}_1}<0$$
	It is evident from combining all of the inequalities that $\tau_1>0$, resulting in the stability of $(0,\hat{v})$. 
\end{proof}

\section{Non-Existence of Coexistence Steady State}\label{section_6}
Let $(u^\ast,v^\ast)$ represent the co-existence steady state. Consequently, the following stationary problem is satisfied by $u^\ast,v^\ast(>0)$.
\begin{align} 
	\label{co-existence_system22}
	\begin{cases}
		\vspace{0.2cm}
		& ({\mathsf{d}_1}{u_x^\ast}-{\alpha_1}{u^\ast})_x+rr_1u^\ast{\left[1-\displaystyle\frac{u^\ast+v^\ast}{K_1}\right]}=0,\quad\displaystyle a<x<b,\ t>0 \\
		\vspace{0.2cm}
		&({\mathsf{d}_2}{v_x^\ast}-{\alpha_2}{v^\ast})_x+rr_1v^\ast{\left[1-\displaystyle\frac{u^\ast+v^\ast}{K_1}\right]}=0,\quad\displaystyle a<x<b,\ t>0 \\
		\vspace{0.2cm}
		&{\mathsf{d}_1}{u_x^\ast}(a,t)-{\alpha_1}u^\ast(a,t)=	{\mathsf{d}_1}{u_x^\ast}(b,t)-{\alpha_1}u^\ast(b,t)=0, \quad\displaystyle t>0,\\
		\vspace{0.2cm}
		&{\mathsf{d}_2}{v_x^\ast}(a,t)-{\alpha_2}v^\ast(a,t)={\mathsf{d}_2}{v_x^\ast}(b,t)-{\alpha_2}v^\ast(b,t)=0,\quad\displaystyle t>0,\\
	\end{cases}
\end{align}

\begin{Lemma} \cite{tisha2}
	\label{lm10}
	For system \eqref{co-existence_system22}, assume that $K_1\not\equiv\text{const},\ \mu\in [0,1),\ \mathsf{d}_1,\mathsf{d}_2>0,\ \alpha_1,\alpha_2\in \mathbb{R}$\, and $(u^\ast,v^\ast)$\ are co-existence steady states. Consequently, for any pair of points $a\leq a_1\leq b_1\leq b$, we acquire
	\begin{multline}
		\label{co-int_12}
		\displaystyle \frac{1}{\mathsf{d}_1}\int_{a_1}^{b_1}{[\mathsf{d}_1-\mathsf{d}_2]\cdot\left[ v_x-\displaystyle\frac{\alpha_1-\alpha_2}{\mathsf{d}_1-\mathsf{d}_2}v\right]\cdot[\mathsf{d}_1u_x-\alpha_1 u]\cdot e^{-\frac{\alpha_1}{\mathsf{d}_1}x}}dx\\
		=\left[(\mathsf{d}_1u_x-\alpha_1 u)\cdot e^{-\frac{\alpha_1}{\mathsf{d}_1}x}\cdot v\right]\displaystyle\bigg\vert_{a_1}^{b_1}-\left[(\mathsf{d}_2 v_x-\alpha_2 v)\cdot e^{-\frac{\alpha_1}{\mathsf{d}_1}x}\cdot u \right]\displaystyle\bigg\vert_{a_1}^{b_1}
	\end{multline}
	and 
	\begin{multline}
		\label{co-int_22}
		\displaystyle \frac{1}{\mathsf{d}_2}\displaystyle\int_{a_1}^{b_1}{[\mathsf{d}_1-\mathsf{d}_2]\cdot\left[ u_x-\displaystyle\frac{\alpha_2-\alpha_1}{\mathsf{d}_2-\mathsf{d}_1}u\right]\cdot[\mathsf{d}_2v_x-\alpha_2 v]\cdot e^{-\frac{\alpha_2}{\mathsf{d}_2}x}}dx\\
		=\left[(\mathsf{d}_2 v_x-\alpha_2 v)\cdot e^{-\frac{\alpha_2}{\mathsf{d}_2}x}\cdot u\right]\displaystyle\bigg\vert_{a_1}^{b_1}-\left[(\mathsf{d}_1 u_x-\alpha_1 u)\cdot e^{-\frac{\alpha_2}{\mathsf{d}_2}x}\cdot v \right]\displaystyle\bigg\vert_{a_1}^{b_1}
	\end{multline}
\end{Lemma}                        
\noindent\textbf{Note:} For simplicity, $(u^\ast,v^\ast)$ is replaced by $(u,v)$. 
\begin{proof}
Rewrite the system's second equation, \eqref{co-existence_system22}, to demonstrate the identity \eqref{co-int_12} as
	\begin{align}
		&	{\mathsf{d}_2}{\Delta v}-{\alpha_2}{ v_x}+rr_1v{\left[1-\displaystyle\frac{u+v}{K_1}\right]}=0\nonumber\\
		\Rightarrow& {\mathsf{d}_2}{\Delta v}-{\alpha_2}{v_x}+{\mathsf{d}_1}{\Delta v}-{\alpha_1}{v_x}-{\mathsf{d}_1}{\Delta v}+{\alpha_1}{ v_x}+rr_1v{\left[1-\displaystyle\frac{u+v}{K_1}\right]}=0\nonumber\\
		\Rightarrow& {\mathsf{d}_1}{\Delta v}-{\alpha_1}{v_x}+rr_1v{\left[1-\displaystyle\frac{u+v}{K_1}\right]}=[({\mathsf{d}_1}-{\mathsf{d}_2}){v_x}+({\alpha_2}-{\alpha_1})]_x{v_x}\nonumber
	\end{align}
	Then the system \eqref{co-existence_system22} becomes
	\begin{align} 
		\begin{cases}
			\vspace{0.2cm}
			& ({\mathsf{d}_1}{u_x}-{\alpha_1}{u})_x+rr_1u{\left[1-\frac{u+v}{K_1}\right]}=0,\quad\displaystyle a<x<b,\ t>0 \\
			\vspace{0.2cm}
			&({\mathsf{d}_1}{v_x}-{\alpha_1}{v})_x+rr_1v{\left[1-\frac{u+v}{K_1}\right]}={(\mathsf{d}_1-\mathsf{d}_2)}_x{v_x}-{(\alpha_2-\alpha_1)}{v_x},\quad\displaystyle a<x<b,\ t>0 
		\end{cases}
	\end{align}
	After multiplying $e^{-\frac{\alpha_1}{\mathsf{d}_1}x}v$ by the first equation, we integrate over $[a,b]$ to obtain
	\begin{align*}
		\displaystyle \int_{a_1}^{b_1}{
			[\mathsf{d}_1 u_x-\alpha_1 u]_x\cdot e^{-\frac{\alpha_1}{\mathsf{d}_1}x}v+rr_1u\left[1-\displaystyle\frac{u+v}{K_1}\right]e^{-\frac{\alpha_1}{\mathsf{d}_1}x}v}dx=0
	\end{align*}
	\vspace{-.3cm}
	\begin{multline}
		\label{co-eq11}
		\Rightarrow [\mathsf{d}_1 u_x-\alpha_1 u]\cdot e^{-\frac{\alpha_1}{\mathsf{d}_1}x}v\displaystyle\bigg \vert_{a_1}^{b_1}-\displaystyle\frac{1}{\mathsf{d}_1}\displaystyle \int_{a_1}^{b_1}{[\mathsf{d}_1 u_x-\alpha_1 u]\cdot [\mathsf{d}_1v_x-\alpha_1 v]\cdot e^{-\frac{\alpha_1}{\mathsf{d}_1}x}}dx+\\
		\displaystyle	\int_{a_1}^{b_1}{rr_1u\left[1-\displaystyle\frac{u+v}{K_1}\right]e^{-\frac{\alpha_1}{\mathsf{d}_1}x}v}dx=0
	\end{multline}
	Multiply the first equation by $e^{-\frac{\alpha_1}{\mathsf{d}_1}x}u$\ and then integrate over $[a_1,b_1]$, we get
	\begin{multline*}
		\displaystyle \int_{a_1}^{b_1}{
			[\mathsf{d}_1v_x-\alpha_1 v]_x\cdot e^{-\frac{\alpha_1}{\mathsf{d}_1}x}u+rr_1v\left[1-\displaystyle\frac{u+v}{K_1}\right]e^{-\frac{\alpha_1}{\mathsf{d}_1}x}u}dx\\
		=\displaystyle \int_{a_1}^{b_1}{\left[(\mathsf{d}_1-\mathsf{d}_2)v_x+(\alpha_2-\alpha_1) v\right]_x e^{-\frac{\alpha_1}{\mathsf{d}_1}x}u}dx
	\end{multline*}\vspace{-.3cm}
	\begin{multline}
		\label{co-eq22}
		\Rightarrow [\mathsf{d}_1v_x-\alpha_1 v]\cdot e^{-\frac{\alpha_1}{\mathsf{d}_1}x}u\displaystyle\bigg \vert_{a_1}^{b_1}-\displaystyle \int_{a_1}^{b_1}{ [\mathsf{d}_1 u_x-\alpha_1 u]\cdot [\mathsf{d}_1v_x-\alpha_1 v]\cdot e^{-\frac{\alpha_1}{\mathsf{d}_1}x}}dx+\\
		\displaystyle\int_{a_1}^{b_1}{rr_1v\left[1-\displaystyle\frac{u+v}{K_1}\right]e^{-\frac{\alpha_1}{\mathsf{d}_1}x}u}dx
		=\left[(\mathsf{d}_1-\mathsf{d}_2)v_x+(\alpha_2-\alpha_1) v\right] e^{-\frac{\alpha_1}{\mathsf{d}_1}x}u\bigg \vert_{a_1}^{b_1}\\
		-\displaystyle\frac{1}{\mathsf{d}_1}\displaystyle\int_{a_1}^{b_1}{[\mathsf{d}_1 u_x-\alpha_1 u]\cdot[\mathsf{d}_1-\mathsf{d}_2]\cdot\left[v_x-\displaystyle\frac{\alpha_1-\alpha_2}{\mathsf{d}_1-\mathsf{d}_2} v\right] \left( e^{-\frac{\alpha_1}{\mathsf{d}_1}x}\right)}dx
	\end{multline}
The identity \eqref{co-int_12} is obtained by subtracting equation \eqref{co-eq22} from equation \eqref{co-eq11}. Similarly, the identity \eqref{co-int_22} can be derived.
\end{proof}

Assume that there is a coexistence solution. Examine the following transformation for the co-existence solution $(u^\ast,v^\ast)$: $$\mathcal{T}:=\displaystyle\frac{u_x}{u}\ \text{and}\ \mathcal{S}:=\displaystyle\frac{v_x}{v}$$
 We are able to obtain \cite{2016}
\begin{align}
	\label{co-system22}
	\begin{cases}
		&-\mathsf{d}_1\Delta \mathcal{T}+[\alpha_1-2\mathsf{d}_1\mathcal{T}]\mathcal{T}_x+u\mathcal{T}+v\mathcal{S}=0, \quad a<x<b,\\
		&-\mathsf{d}_2\Delta \mathcal{S}+[\alpha_2-2\mathsf{d}_1\mathcal{S}]\mathcal{S}_x+u\mathcal{T}+v\mathcal{S}=0, \quad a<x<b,\\
		&\mathcal{T}(a)= \mathcal{T}(b)=\frac{\alpha_1}{\mathsf{d}_1}>0,\\
		&\mathcal{S}(a)=\mathcal{S}(b)=\frac{\alpha_2}{\mathsf{d}_2}>0,
	\end{cases}
\end{align}
Recall the system \eqref{harvest_free}
\begin{align} 
	\begin{cases}
		\vspace{0.2cm}
		&-({\mathsf{d}_1}{u_x}-{\alpha_1}{u})_x=rr_1u{\left[1-\displaystyle\frac{u+v}{K_1}\right]},\quad\displaystyle a<x<b,\ t>0 \\
		\vspace{0.2cm}
		&-({\mathsf{d}_2}{v_x}-{\alpha_2}{v})_x=rr_1v{\left[1-\displaystyle\frac{u+v}{K_1}\right]},\quad\displaystyle a<x<b,\ t>0 
	\end{cases}
\end{align}
which can be expressed as
\begin{align} 
	\begin{cases}
		\vspace{0.2cm}
		&\displaystyle\frac{-({\mathsf{d}_1}{\Delta u}-{\alpha_1}{u_x})}{u}=rr_1{\left[1-\displaystyle\frac{u+v}{K_1}\right]}, \\
		\vspace{0.2cm}
		&\displaystyle\frac{-({\mathsf{d}_2}{\Delta v}-{\alpha_2}{v_x})}{v}=rr_1{\left[1-\displaystyle\frac{u+v}{K_1}\right]}, 
	\end{cases}
\end{align}
This implies
\begin{align*}
	\vspace{0.2cm}
	&\displaystyle\frac{-({\mathsf{d}_1}{\Delta u}-{\alpha_1}{u_x})}{u}=\displaystyle\frac{-({\mathsf{d}_2}{\Delta v}-{\alpha_2}{v_x})}{v}=rr_1{\left[1-\displaystyle\frac{u+v}{K_1}\right]}\\
	\vspace{0.2cm}
	\Rightarrow& -\mathsf{d}_1\frac{\Delta u-{u_x}^2}{u^2}+\alpha_1 \displaystyle\frac{u_x}{u}-\mathsf{d}_1\displaystyle\frac{{u_x}^2}{u^2}=-\mathsf{d}_2\displaystyle\frac{\Delta v-{v_x}^2}{v^2}+\alpha_2 \displaystyle\frac{v_x}{v}-\mathsf{d}_2\displaystyle\frac{{v_x}^2}{v^2}
	\vspace{0.2cm}
\end{align*}
Therefore
\begin{equation}
	\label{non co-ex22}
	-\mathsf{d}_1 \Delta \mathcal{T}+\alpha_1 \mathcal{T}-\mathsf{d}_1 \mathcal{T}^2=-\mathsf{d}_2 \Delta \mathcal{S}+\alpha_2 \mathcal{S}-\mathsf{d}_2 \mathcal{S}^2, \qquad a<x<b.
\end{equation}

\begin{Lemma}
	\label{lm27}
	Let $\mathcal{S}$\ and $\mathcal{T}$\ be described as previously, and the identity \eqref{non co-ex22} is accurate. Therefore, the following propositions cannot be true for any two points $a\leq a_1\leq b_1\leq b$.
	\begin{enumerate}
		\item[(a)] The local maximum of $\mathcal{T}$\ is positive in $(a_1,b_1)$\, and $\mathcal{S}\geq 0$\ in $[a_1,b_1]$;
		\item[(a)] The local maximum of $\mathcal{S}$\ is positive in $(a_1,b_1)$\, and $\mathcal{T}\geq 0$\ in $[a_1,b_1]$;
		\item[(c)] The local minimum of $\mathcal{T}$\ is negative in $(a_1,b_1)$\, and $\mathcal{S}\leq 0$\ in $[a_1,b_1]$;
		\item[(d)] The local minimum of $\mathcal{S}$\ is negative in $(a_1,b_1)$\, and $\mathcal{T}\leq 0$\ in $[a_1,b_1]$;
	\end{enumerate}
\end{Lemma}
\begin{proof}
Assuming $p\in (a_1,b_1)$ and that $\mathcal{T}$\ has a positive local maximum in $(a_1,b_1)$\, then $$\mathcal{T}(p) >0,\ \mathcal{T}^\prime (p)=0,\ \text{and}\ \mathcal{T}^{\prime\prime}\leq 0.$$This suggests that $\mathcal{S}(p)<0$\ at \ from the system's first equation \eqref{co-system22} at the point $p$, which runs counter to our lemma's assertion.\\
Similarly, it is possible to prove the other claims.
\end{proof}
\noindent Using the earlier findings, it is now possible to prove the non-existence of coexistence stable state \cite{tisha2,he3,zhou-2016, tisha1}.
\begin{Lemma}
	\label{lm26}
Assume that \eqref{c_1} holds and $K_1\not\equiv \text{const},\ \mu\in [0,1),\ \mathsf{d}_1>\mathsf{d}_2>0,\ \alpha_1>\alpha_2>0$. In such case, there is no coexistence stable state in system \eqref{harvest_system2}.
\end{Lemma}
\begin{proof}
	Let us define
	$$\mathcal{A}(x)=\mathsf{d}_1u_x-\alpha_1 u\ \text{and} \ \mathcal{B}(x)=\mathsf{d}_2 v_x-\alpha_2v$$
	By using boundary condition\ $\mathcal{A}(a)=\mathcal{A}(b)=0\ \text{and}\ \mathcal{B}(a)=\mathcal{B}(b)=0$. Also
	$$\mathcal{A}^\prime(x)=-rr_1u\left[1-\displaystyle\frac{u+v}{K_1}\right]\ \text{and}\ \mathcal{B}^\prime(x)=-rr_1v\left[1-\displaystyle\frac{u+v}{r_1K}\right]$$
	\textbf{Step 1:} $\mathcal{A}^\prime(x)>0\Leftrightarrow \mathcal{B}^\prime(x)>0,\ \mathcal{A}^\prime(x)<0\Leftrightarrow\mathcal{B}^\prime(x)<0,\ \mathcal{A}^\prime(x)=0\Leftrightarrow \mathcal{B}^\prime(x)=0$\ in $[a,b]$.\\[2ex]
	It is evident that the same function $\left[1-\displaystyle\frac{u+v}{K_1}\right]$ determines the sign of $\mathcal{A}^\prime(x)$\ and \ $\mathcal{B}^\prime(x)$ given that $u,v$\, and $r$\ are all strictly positive. Step 1 is thus fulfilled.\\[2ex]
	\textbf{Step 2:} In any interval $[a_1,b_1]\subset[a,b]$, $\mathcal{A}$\ and $\mathcal{B}$\ cannot be identically zero if $\mathsf{d}_1=\mathsf{d}_2=\mathsf{d}$. Pretend that $\mathcal{A}\equiv 0$ in any interval $[a_1,b_1]\subset[a,b]$, if at all possible. Next,
	\begin{align*}
		&\mathsf{d}u_x-\alpha_1 u\equiv 0\\
		&\Rightarrow u_x\equiv \displaystyle\frac{\alpha_1}{\mathsf{d}}u \quad \text{in}\ [a_1,b_1]\\
		& \Rightarrow u(x)=u(a_1)e^{\frac{\alpha_1}{\mathsf{d}}}\quad \text{for}\ x\in[a_1,b_1]
	\end{align*}
	Conversely, because $\mathcal{A}\equiv 0$\ in\ $[a_1,b_1]\subset[a,b],\ \mathcal{A}^\prime(x) \equiv 0 \ \text{for}\ x\in[a_1,b_1]$. Then Step 1 implies that $1-\displaystyle\frac{u+v}{K_1}\equiv 0\ \text{in}\ [a_1,b_1]$,\ so \ $\mathcal{B}^\prime(x) \equiv 0\Leftrightarrow \mathcal{B}\triangleq \mathsf{d}_2v_x-\alpha_2v\equiv \text{const}  $. Now $v(x)=r-u(x)=r-u(a_1)e^{\frac{\alpha_1}{\mathsf{d}}}$. Then $\mathcal{B}\equiv \text{const}\Leftrightarrow -\mathsf{d}u_x-\alpha_2\left(r-u(a_1)e^{\frac{\alpha_1}{\mathsf{d}}}\right)\equiv \text{const}\Leftrightarrow -\mathsf{d}\displaystyle\frac{\alpha_1}{\mathsf{d}}u-\alpha_2\left(r-u(a_1)e^{\frac{\alpha_1}{\mathsf{d}}x}\right)\equiv \text{const}\Leftrightarrow -\alpha_1u_x-0+\alpha_2\displaystyle\frac{\alpha_1}{\mathsf{d}}u(a_1)e^{\frac{\alpha_1}{\mathsf{d}}x}=0\Leftrightarrow  -\alpha_1u_x+\displaystyle\frac{\alpha_1}{\mathsf{d}}\alpha_2 u=0\Leftrightarrow \alpha_1u_x=\alpha_2u_x\Leftrightarrow \alpha_1=\alpha_2$, which is a contradiction of the assumption of our lemma.\\[2ex]
	\textbf{Step 3:}
In $(a,\delta)\cup (b-\delta,b)$, $\mathcal{A}<0$\ and $\mathcal{B}<0$\ in\ $(a,\delta)\cup (b-\delta,b)$\ if a small $\delta>0$ exists.\\[2ex]Since $\mathcal{A}(a)=\mathcal{A}(b)=0$, if $\mathcal{A}<0$\ in\ $(a,\delta)$\ and\ $(b-\delta,b)$, then if required, one can reduce $\delta>0$ to find that\ $\mathcal{A}$\ decreases in $(a,\delta)$ and increases in \ $(b-\delta,b)$. Additionally, since $\mathcal{B}(a)=\mathcal{B}(b)=0$, Step 1 suggests that $\mathcal{B}<0$\ in\ $(a,\delta)$\ and\ $(b-\delta,b)$. We can work with $(a,\delta)$\ just for convenience because $\mathcal{A}$\ behaves similarly in\ $(a,\delta)$\ and\ $(b-\delta,b)$. This suggests that $$\mathcal{A}<0\ \text{in}\ (a,\delta)\ \text{for \ small}\ \delta>0,$$ If required, one can then reduce $\delta>0$ to determine from Step 1 that $$\mathcal{B}<0\ \text{in} \ (a,\delta)\ \text{for\ small}\ \delta>0$$ Let the first zero point of $\mathcal{A}$\ and $\mathcal{B}$\ in $[\delta,b]$\ be\ $p$\ and\ $q$\, respectively ($\mathcal{A}(a)=\mathcal{B}(b)=0$\ ensures the existence of zero points). Assume that $p\leq q$.
	$$\mathcal{A}(a)=\mathcal{A}(p)=0,\ \mathcal{A}(x)>0, \quad x\in(a,p)$$
	Because $p\leq q$, 
	$$\mathcal{B}(a)=0\leq\mathcal{B}(p),\ \mathcal{B}(x)>0, \quad x\in(a,p)$$
	So that we have
	\begin{align}
		\label{heq11}
		\mathcal{T}(0)=\mathcal{T}(p)=\displaystyle\frac{\alpha_1}{\mathsf{d}_1}, \qquad \mathcal{T}(x)>\displaystyle\frac{\alpha_1}{\mathsf{d}_1}, \quad x\in(a,p)
	\end{align}
	and 
	\begin{align}
		\label{heq2}
		\mathcal{S}(0)=\displaystyle\frac{\alpha_2}{\mathsf{d}_2}\leq \mathcal{S}(p), \qquad \mathcal{S}(x)>\displaystyle\frac{\alpha_2}{\mathsf{d}_2}>0, \quad x\in(a,p)
	\end{align}
	with $\mathcal{S}>0$, demonstrating that $\mathcal{T}$\ admits a positive local maximum in $(a,p)$. But Lemma \ref{lm27}'s statement(a) states that this cannot happen. Thus, Step 3 is fulfilled.\\[2ex]
	\textbf{Step 4:} In $(a,b)$, both $\mathcal{A}$\ and $\mathcal{B}$\ must change sign.\\[2ex]
	Let $p$ be $\mathcal{B}$'s initial zero point in $(a,b)$ (the existence of $p$ is guaranteed by $\mathcal{B}(b)=0)$. Because of Step 3 
	$$\mathcal{B}(x)=\mathsf{d}_2v_x-\alpha_2v<0\Leftrightarrow v_x-\displaystyle\frac{\alpha_2}{\mathsf{d}_2}v<0,\quad x\in(a,p)$$
As stated in Lemma \ref{lma4},
	$$\displaystyle\frac{\alpha_1-\alpha_2}{\mathsf{d}_1-\mathsf{d}_2}\geq\displaystyle\frac{\alpha_2}{\mathsf{d}_2}$$
	so
	$$v_x-\displaystyle\frac{\alpha_1-\alpha_2}{\mathsf{d}_1-\mathsf{d}_2}v<0,\quad x\in(a,p)$$
	Suppose for contradiction that
	$$\mathcal{A}=\mathsf{d}_1u_x-\alpha_1u\leq 0,\quad x\in[a,p]$$
	Configuring $a_1=a,b_1=p$\ in \eqref{co-int_12}\ we obtain
	\begin{align*}
		&0<	\displaystyle \frac{1}{\mathsf{d}_1}\int_{a}^{p}{[\mathsf{d}_1-\mathsf{d}_2]\cdot\left[v_x-\displaystyle\frac{\alpha_1-\alpha_2}{\mathsf{d}_1-\mathsf{d}_2}v\right]\cdot[\mathsf{d}_1u_x-\alpha_1 u]\cdot e^{-\frac{\alpha_1}{\mathsf{d}_1}x}}dx\\
		&	=\left[(\mathsf{d}_1u_x-\alpha_1 u)\cdot e^{-\frac{\alpha_1}{\mathsf{d}_1}x}\cdot v\right]\displaystyle\bigg\vert_{p}-\left[(\mathsf{d}_2v_x-\alpha_2 v)\cdot e^{-\frac{\alpha_1}{\mathsf{d}_1}x}\cdot u \right]\displaystyle\bigg\vert_{p}\leq 0
	\end{align*}
In addition to demonstrating that $\mathcal{A}$ must change sign in $(a,p)$, this contradiction also demonstrates that $\mathcal{A}$\ must change sign before $\mathcal{B}$. This step must be completed by confirming that $\mathcal{B}$ likewise changes sign in $(a,b)$. To accomplish this, let q be $\mathcal{A}$'s final zero point. In line with Step 3,
	$$\mathcal{A}=\mathsf{d}_1u_x-\alpha_1u\leq 0,\quad x\in(q,b)$$
	As stated in Lemma \ref{lma4}, 
	$$\displaystyle\frac{\alpha_1-\alpha_2}{\mathsf{d}_1-\mathsf{d}_2}\geq\displaystyle\frac{\alpha_1}{\mathsf{d}_1}$$
	and so
	$$u_x-\frac{\alpha_1-\alpha_2}{\mathsf{d}_1-\mathsf{d}_2}u<0,\quad x\in(q,b)$$
	If possible $$\mathcal{B}=\mathsf{d}_2v_x-\alpha_2v\leq 0,\quad x\in(q,b)$$
	Setting $a_1=q,~ b_1=b$\ in \eqref{co-int_22}\ we obtain
	\begin{align*}
		&0>\displaystyle \frac{1}{\mathsf{d}_2}\int_{a_1}^{b_1}{[\mathsf{d}_1-\mathsf{d}_2]\cdot\left[u_x-\displaystyle\frac{\alpha_2-\alpha_1}{\mathsf{d}_2-\mathsf{d}_1}u\right]\cdot[\mathsf{d}_2v_x-\alpha_2 v]\cdot e^{-\frac{\alpha_2}{\mathsf{d}_2}x}}dx\\
		&=\left[(\mathsf{d}_2v_x-\alpha_2 v)\cdot e^{-\frac{\alpha_2}{\mathsf{d}_2}x}\cdot u\right]\displaystyle\bigg\vert_q-\left[(\mathsf{d}_1u_x-\alpha_1 u)\cdot e^{-\frac{\alpha_2}{\mathsf{d}_2}x}\cdot v \right]\displaystyle\bigg\vert_q\geq 0
	\end{align*}
which contradicts itself. This indicates that $(q,b)$ requires\ $\mathcal{B}$\ to change sign. Step 4 is therefore valid. $\mathcal{A}$\ has at least two zero points in $(a,b)$, according to the instructions above. If $$a=a_0<a_1<a_2<\cdots <a_n=b, \qquad n\geq 3(n\in \mathbb{Z}),$$ represents these zero points, then it is evident that $\mathcal{B}<0$\ in $(a_0,a_1]$\ because of Step 4.
	Step 3 indicates that there is\ $a<a_1<a_2<b$\ such that 
	$$\mathcal{B}(a)=\mathcal{B}(a_1)=\mathcal{B}(a_2)=0, \quad \mathcal{B}(x)\leq 0 \ \text{in}\ (a,a_1), \quad \mathcal{B}(x)> 0 \ \text{in}\ (a_1,a_2)$$
	and so
	$$\mathcal{S}(x)\geq \displaystyle\frac{\alpha_2}{\mathsf{d}_2}>0\ \text{in}\ [a_1,a_2]$$\\[2ex]
	\textbf{Step 5:} In\ $(a,a_1]$, $\mathcal{A}$\ is strictly negative. If not, there is at least one zero point in\ $(a,a_1]$ of\ $\mathcal{A}$. Assume that the first zero point of $\mathcal{A}$ is $a_0$. Next,
	\begin{align}
		\label{eq11}
		\mathcal{A}(a)=\mathcal{A}(a_0)=0,~ \mathcal{A}(x)<0 \ \text{in} \ (a,a_0)\subset (a,a_1] 
	\end{align}
	which indicates that 
	\begin{align}
		\label{eq22}
		&\displaystyle \frac{1}{\mathsf{d}_1}[\mathsf{d}_1u_x-\alpha_1 u]\cdot e^{-\frac{\alpha_1}{\mathsf{d}_1}x}<0
	\end{align}
	It is clear that,\ $\mathcal{B}(a_0)<0$\ and 
	\begin{align}
		\label{eq33}
		&	\mathsf{d}_1v_x(a_0)<\alpha_1 v(a_0)
	\end{align}
	Also by Step 3, \ $\mathcal{B}(x)>0$\ in \ $(a,a_0)\subset (a,a_1] $. And so 
	\begin{align}
		\label{eq44}
		&v_x-\displaystyle\frac{\alpha_1-\alpha_2}{\mathsf{d}_1-\mathsf{d}_2}v>0,\quad x\in(a, a_0)
	\end{align}
	Setting $a_1=a,b_1=a_0$\ in \eqref{co-int_12}, we obtain
	\begin{align}
		\label{eq55}
		&	\displaystyle \frac{1}{\mathsf{d}_1}\displaystyle\int_{a}^{a_0}{[\mathsf{d}_1-\mathsf{d}_2]\cdot\left[v_x-\displaystyle\frac{\alpha_1-\alpha_2}{\mathsf{d}_1-\mathsf{d}_2}v\right]\cdot[\mathsf{d}_1u_x-\alpha_1 u]\cdot e^{-\frac{\alpha_1}{\mathsf{d}_1}x}}dx\nonumber\\
		&	=\left[(\mathsf{d}_1u_x-\alpha_1 u)\cdot e^{-\frac{\alpha_1}{\mathsf{d}_1}x}\cdot v\right]\displaystyle\bigg\vert_{a}^{a_0}-\left[(\mathsf{d}_2v_x-\alpha_2 v)\cdot e^{-\frac{\alpha_1}{\mathsf{d}_1}x}\cdot u \right]\displaystyle\bigg\vert_{a}^{a_0}
	\end{align}
	Applying \eqref{eq11}, it becomes 
	\begin{align}
		&	\displaystyle \frac{1}{\mathsf{d}_1}\displaystyle\int_{a}^{a_0}{[\mathsf{d}_1-\mathsf{d}_2]\cdot\left[v_x-\displaystyle\frac{\alpha_1-\alpha_2}{\mathsf{d}_1-\mathsf{d}_2}v\right]\cdot[\mathsf{d}_1u_x-\alpha_1 u]\cdot e^{-\frac{\alpha_1}{\mathsf{d}_1}x}}dx\nonumber\\
		&	=-\left[(\mathsf{d}_2v_x-\alpha_2 v)\cdot e^{-\frac{\alpha_1}{\mathsf{d}_1}x}\cdot u \right]\displaystyle\bigg\vert_{a}^{a_0}
	\end{align}
	Clearly,
	\begin{align}
		&	\displaystyle \frac{1}{\mathsf{d}_1}\displaystyle\int_{a}^{a_0}{[\mathsf{d}_1-\mathsf{d}_2]\cdot\left[v_x-\displaystyle\frac{\alpha_1-\alpha_2}{\mathsf{d}_1-\mathsf{d}_2}v\right]\cdot[\mathsf{d}_1u_x-\alpha_1 u]\cdot e^{-\frac{\alpha_1}{\mathsf{d}_1}x}}dx\geq 0
	\end{align}
	Because of inequality \eqref{eq22}, this is contradicted by inequality \eqref{eq33} and \eqref{eq44}. Therefore, Step 5 is valid.\\[2ex]
	\textbf{Step 6:} \ In $(a_1,a_2)$, $\mathcal{B}$\ is non positive, whereas $\mathcal{B}(a_2)$\ is strictly negative and $\mathcal{A}$\ is strictly positive.\\[2ex]
	By \eqref{heq11}, 
	\begin{align}
		& \mathcal{T}(a_1)=\mathcal{T}(a_2)=\displaystyle\frac{\alpha_1}{\mathsf{d}_1},\qquad \mathcal{T}>\displaystyle\frac{\alpha_1}{\mathsf{d}_1}\ \text{in}\ (a_1,a_2)
	\end{align}
	Let us suppose that a positive local maximum in\ $ (a_1,a_2)$ must be admitted by $\mathcal{S}$. We examine the following two situations to support this notion.
	$$\text{Case (i):}\ \mathcal{B}\leq 0$$
	$$\text{Case (ii):}\ \mathcal{B}>0$$
	For Case (i), there exists two numbers \ $p,q$\ and \ $a_1<p<q<a_2$\ such that
	$$\mathcal{B}(p)=\mathcal{B}(q)=0,\quad \mathcal{B}>0 \ \text{in}\ (p,q)$$
Likewise,
	$$\mathcal{S}(p)=\mathcal{S}(q)=\displaystyle\frac{\alpha_2}{\mathsf{d}_2},\quad \mathcal{S}>\displaystyle\frac{\alpha_2}{\mathsf{d}_2} \ \text{in}\ (p,q)$$
this suggests our presumption. There is \ $a_1<s<a_2$\ for Case (ii) such that
	$$\mathcal{B}(a_1)=\mathcal{B}(s)=0,\quad \mathcal{B}>0 \ \text{in}\ (s,a_2]$$
	equivalently,
	\begin{align}
		\label{heq33}
		\mathcal{S}(s)=\displaystyle\frac{\alpha_2}{\mathsf{d}_2},\quad \mathcal{S}>\displaystyle\frac{\alpha_2}{\mathsf{d}_2} \ \text{in}\ (s,a_2]
	\end{align}
	Now since \ $ \mathcal{A}> 0$\ in \ $(a_1,a_2)$\ and \ $ \mathcal{A}(a_2)= 0$, thus, at $x\rightarrow a_2^-$, $ \mathcal{A}$\ exhibits a declining behaviour. Specifically, $\mathcal{A}(a_2^-)\leq 0$. Thus by step 1,\ $\mathcal{B}(a_2^-)\leq 0$.
	Notice that 
	\begin{align}
		& \mathcal{B}^\prime(a_2^-)=\mathsf{d}_2\Delta v-\alpha_2 v_x\displaystyle\big\vert_{a_2}=\displaystyle \mathsf{d}_2\frac{v\Delta v-\displaystyle\frac{\alpha_2}{\mathsf{d}_2}vv_x}{v}\displaystyle\bigg\vert_{a_2},
	\end{align}
	\begin{align}
		&	\mathcal{S}^\prime(a_2^-)=\left[\displaystyle\frac{v_x}{v}\right]_x=\displaystyle\frac{v\Delta v-{v_x}^2}{v^2}.
	\end{align}
	Since $\displaystyle\frac{v_x}{v}>\displaystyle\frac{\alpha_2}{\mathsf{d}_2}>0$\ in\ $(a_1,a_2)\ (\mathcal{B}>0)$ by \eqref{co-system22}, so at\ $x\rightarrow a_2^-$\
	\begin{align}
		& v\Delta v-{v_x}^2\displaystyle\big\vert_{a_2}<v\Delta v-\displaystyle\frac{\alpha_2}{\mathsf{d}_2}vv_x\displaystyle\big\vert_{a_2}=\displaystyle\frac{\mathcal{B}^\prime(a_2)v(a_2)}{\mathsf{d}_2}\leq0\ \left(\mathcal{B}^\prime({a_2^-})\leq0\right).
	\end{align}
Together with \eqref{heq33}, we thus obtain $$\mathcal{S}^\prime(a_2^-)\leq 0$$ which suggests that $\mathcal{S}$ must permit a positive local maximum in $(s,a_2)$. Thus, our hypothesis is correct. However, this contradicts Lemma \ref{lm27}'s statement (a). Consequently, since \ $\mathcal{A}(a_2)<0$, by Step 1,\ $\mathcal{B}(a_2)<0$, $$\mathcal{B}\leq 0\ \text{in}\ (a_1,a_2)$$ Thus, Step 6 is validated.\\[2ex]
	\textbf{Step 7:}\ In $(a_1,a_2)$\, $\mathcal{B}$\ is not positive, whereas $\mathcal{B}(a_2)$\ is strictly negative, and $\mathcal{A}$\ is strictly negative in\ $(a_1,a_2)$.\\[2ex]
	The fact that $\mathcal{B}$\ is strictly negative in $(a_1,a_2]$ is readily demonstrable. If not, there must be zero points in\ $\mathcal{B}$. Assume that $r$ represents one of the points that is closest to $a_1$. Then
	\begin{align}
		\mathcal{A}(a_1)=0,~ \mathcal{A}(r)\leq 0,\ \mathcal{A}<0\ \text{in}\ (a_1,r), 
	\end{align}
	and so
	\begin{align}
		\mathcal{B}(a_1)<0,~ \mathcal{B}(r)= 0,\ \mathcal{B}<0 \ \text{in}\ (a_1,r), 
	\end{align}
By Step 4, this suggests
	\begin{align}
		&v_x-\displaystyle\frac{\alpha_1-\alpha_2}{\mathsf{d}_1-\mathsf{d}_2}v<0,\quad \text{in} \ (a_1,r),
	\end{align}
	Setting $\alpha_1 =a_1,~b_1=r$\ in \eqref{co-int_12}\ we obtain
	\begin{align*}
		&0<	\displaystyle \frac{1}{\mathsf{d}_1}\displaystyle\int_{a_1}^{r}{[\mathsf{d}_1-\mathsf{d}_2]\cdot\left[v_x-\displaystyle\frac{\alpha_1-\alpha_2}{\mathsf{d}_1-\mathsf{d}_2}v\right]\cdot[\mathsf{d}_1u_x-\alpha_1 u]\cdot e^{-\frac{\alpha_1}{\mathsf{d}_1}x}}dx\\
		&	=\left[(\mathsf{d}_1u_x-\alpha_1 u)\cdot e^{-\frac{\alpha_1}{\mathsf{d}_1}x}\cdot v\right]\displaystyle\bigg\vert_{r}-\left[(\mathsf{d}_2v_x-\alpha_2 v)\cdot e^{-\frac{\alpha_1}{\mathsf{d}_1}x}\cdot u \right]\displaystyle\bigg\vert_{r}\\
		& =\mathcal{A}\cdot e^{-\frac{\alpha_1}{\mathsf{d}_1}x}\cdot v\displaystyle\bigg\vert_{r} +\mathcal{B}\cdot e^{-\frac{\alpha_1}{\mathsf{d}_1}x}\cdot u\displaystyle\bigg\vert_{a_1} < 0
	\end{align*}
which contradicts itself. Therefore, Step 7 is valid.\\
This lemma can now be concluded. The next zero point $a_3$\ of $\mathcal{A}$ is taken into consideration once $x$\ passes through $a_2$. We have $$\mathcal{B}\leq 0 \ \text{in}\ (a_2,a_3)\ \text{and}\ \mathcal{B}(a_3)< 0 $$ based on step 4. Similarly, we get $$\mathcal{B}\leq 0 \ \text{in}\ (a_{k-1},a_k), \text{where}\ 1\leq j\leq n,~ j\in \mathbb{Z}$$ So $$\mathcal{B}\leq 0 \ \text{in}\ (a,b)$$ however, this runs counter to Step 4. The conclusion demonstrates that the coexistence stable state $(u^\ast,v^\ast)$ does not exist, proving that our initial assumption was incorrect. The lemma is so proved.
\end{proof}
\begin{Th}
	\label{th_harvest}
	\begin{enumerate}
		\item The global asymptotic stability of $(0,\hat{v})$\ is guaranteed if $K\not\equiv \text{const},\ \mu\in [0,1),\ \mathsf{d}_1>\mathsf{d}_2>0,\ \alpha_1>\alpha_2>0$, and \eqref{c_1} holds.
		\item When $K\not\equiv \text{const},\ \mu\in [0,1],\ \mathsf{d}_1>\mathsf{d}_2>0,\ \alpha_1>\alpha_2>0$ and \eqref{c_2} are true, then for $\omega_1=\displaystyle\frac{\mathsf{d}_2}{\mathsf{d}_1}$, there are two small positive numbers $\epsilon_1$\ and $\epsilon_2$\ such that
		\begin{enumerate}
			\item $(0,\hat{v})$\ is globally asymptotically stable for $\alpha_2$\ $\in$$\left(\omega_1\alpha_1,\omega_1\alpha_1+\epsilon_1\right)$.
			\item $(\hat{u},0)$\ is globally asymptotically stable for $\alpha_2$\ $\in$$\left(\alpha_1-\epsilon_2,\alpha_1\right)$.
			\item The system \eqref{harvest_system2} admits a co-existence steady state for $\alpha_2$ $\in$$\left[\omega_1\alpha_1+\epsilon_1,\alpha_1-\epsilon_2\right]$.
		\end{enumerate}
	\end{enumerate}
\end{Th}
\begin{proof}
	We have demonstrated local stability in Lemma \ref{lm24} and \ref{lm25}, and the non-existence of a co-existence stable state in \ref{lm26}. The first assertion of Theorem \ref{th_harvest} may be obtained simply by applying the theory of competitive systems \cite{tisha2, competitive system1,competitive system2}.
	
To demonstrate the assertion that $2(a)$ First, we note that if \eqref{c_1} holds, that is, if 
	\begin{align*}
		&\displaystyle\frac{\alpha_1}{\mathsf{d}_1}= \displaystyle\frac{\alpha_2}{\mathsf{d}_2} \Rightarrow \alpha_2=\displaystyle\frac{\mathsf{d}_2}{\mathsf{d}_1}\alpha_1
	\end{align*}
This indicates that $(0,\hat{v})$ is globally asymptotically stable.\\[2ex] 
	We will now prove that when\ $\alpha_2\rightarrow\frac{\mathsf{d}_2}{\mathsf{d}_1}\alpha_1^{+}$, there is no coexistence steady state for some small $\epsilon_1>0$\ and \ $\alpha_2\in \left( \displaystyle\frac{\mathsf{d}_2}{\mathsf{d}_1}\alpha_1,\displaystyle\frac{\mathsf{d}_2}{\mathsf{d}_1}\alpha_1+\epsilon_1\right)$. However, the local stability of $(0,\hat{v})$ will constrain because of the constant reliance on the principal eigenvalue on the advection rate of the second species, $\alpha_2$. \\	
When \ $\alpha_2\rightarrow\frac{\mathsf{d}_2}{\mathsf{d}_1}\alpha_1^{+}$, one steady state of coexistence is represented by the symbol $(u_{\alpha_2}^{\ast},v_{\alpha_2}^{\ast})$. For elliptic equations, standard $L_p$ regularity theory \cite{extra} states that if required, subsequences of $\alpha_2$ must be passed to, so that for\ $u^\ast,v^\ast>0$. It is reasonable to presume that
	$$(u_{\alpha_2}^{\ast},v_{\alpha_2}^{\ast})\rightarrow (u^{\ast},v^{\ast}),\ \text{in}\ \ W^{1,p}([a,b])\times W^{1,p}([a,b]) $$
	The Schauder regularity theory and the Sobolev embedding theorem \cite{extra} thus suggest that
	$$(u_{\alpha_2}^{\ast},v_{\alpha_2}^{\ast})\rightarrow (u^{\ast},v^{\ast}),\ \text{in}\ \ C^2([a,b])\times C^2([a,b]) $$
	and fulfill the subsequent set of equations
	\begin{align} 
		\begin{cases}
			\vspace{0.2cm}
			& ({\mathsf{d}_1}{ {u_x^\ast}}-{\alpha_1}{u^\ast})_x+rr_1u^\ast{\left[1-\displaystyle\frac{u^\ast+v^\ast}{K_1}\right]}=0,\quad\displaystyle a<x<b,\ t>0 \\
			\vspace{0.2cm}
			& \frac{\mathsf{d}_2}{\mathsf{d}_1}\left[({\mathsf{d}_1}{ v_x^\ast}-{\alpha_1}{v^\ast})_x\right]+rr_1v^\ast{\left[1-\displaystyle\frac{u^\ast+v^\ast}{K_1}\right]}=0,\quad\displaystyle a<x<b,\ t>0 \\
			\vspace{0.2cm}
			&{\mathsf{d}_1}{u_x^\ast}(a,t)-{\alpha_1}u^\ast(a,t)=	{\mathsf{d}_1}{u_x^\ast}(b,t)-{\alpha_1}u^\ast(b,t)=0, \quad\displaystyle t>0,\\
			\vspace{0.2cm}
			&{\mathsf{d}_1}{v_x^\ast}(a,t)-{\alpha_1}v^\ast(a,t)={\mathsf{d}_1}{v_x^\ast}(b,t)-{\alpha_1}v^\ast(b,t)=0,\quad\displaystyle t>0,
		\end{cases}
	\end{align}
	The existence of $u^\ast,v^\ast>0$\ is denied by Lemma \ref{lm26}. Furthermore, we assert that \ $u^\ast=v^\ast=0$\ is nonexistent. Suppose 
	$$\tilde{u}^{\alpha_2}=\displaystyle \frac{u_{\alpha_2}^{\ast}}{\| u_{\alpha_2}^{\ast}\|_{L^\infty}}\ \text{and}\ \tilde{v}^{\alpha_2}=\displaystyle \frac{v_{\alpha_2}^{\ast}}{\| v_{\alpha_2}^{\ast}\|_{L^\infty}}$$
According to the standard $L_p$ regularity theory \cite{extra} for elliptic equations, we can assume that, if required, we will move to a subsequence of $\alpha_2$
	$$(\tilde{u}_{\alpha_2}^{\ast},\tilde{v}_{\alpha_2}^{\ast})\rightarrow (\tilde{u}^{\ast},\tilde{v}^{\ast}),\ \text{in}\ \ W^{1,p}([a,b])\times W^{1,p}([a,b]) $$
The Schauder regularity theory and the Sobolev embedding theorem \cite{extra} thus suggest that
	$$(\tilde{u}_{\alpha_2}^{\ast},\tilde{v}_{\alpha_2}^{\ast})\rightarrow (\tilde{u}^{\ast},\tilde{v}^{\ast}),\ \text{in}\ \ C^2([a,b])\times C^2([a,b]) $$
and fulfil the subsequent set of equations
	\begin{align} 
		\label{thm1h_int}
		\begin{cases}
			\vspace{0.2cm}
			& ({\mathsf{d}_1}{\tilde{u}_x^\ast}-{\alpha_1}{\tilde{u}^\ast})_x+rr_1\tilde{u}^\ast=0,\quad\displaystyle a<x<b,\ t>0 \\
			\vspace{0.2cm}
			& \frac{\mathsf{d}_2}{\mathsf{d}_1}\displaystyle\left[({\mathsf{d}_1}{\nabla \tilde{v}^\ast}-{\alpha_1}{\tilde{v}^\ast})_x\right]+rr_1\tilde{v}^\ast=0,\quad\displaystyle a<x<b,\ t>0 \\
			\vspace{0.2cm}
			&{\mathsf{d}_1}{\tilde{u}_x^\ast}(a,t)-{\alpha_1}\tilde{u}^\ast(a,t)=	{\mathsf{d}_1}{\tilde{u}_x^\ast}(b,t)-{\alpha_1}\tilde{u}^\ast(b,t)=0, \quad\displaystyle t>0,\\
			\vspace{0.2cm}
			&{\mathsf{d}_1}{\tilde{v}_x^\ast}(a,t)-{\alpha_1}\tilde{v}^\ast(a,t)={\mathsf{d}_1}{\tilde{v}_x^\ast}(b,t)-{\alpha_1}\tilde{v}^\ast(b,t)=0,\quad\displaystyle t>0,\\
		\end{cases}
	\end{align}
	Integrating \eqref{thm1h_int} and imposing boundary condition we get
	$$\displaystyle \int_{a}^{b}{\tilde{u}^\ast rr_1}dx=0$$
	However, we know that $\left\vert \tilde{u}^{\ast}\right\vert_{L^\infty}=1$, that is\ $\tilde{u}^\ast$, in accordance with the norm condition cannot be zero. This results in a contradiction, meaning that our assertion is accurate. This indicates either
	\begin{align}
		\label{suppli11}
		u^\ast=0,\ v^\ast>0 \ \text{in}\ [a,b]
	\end{align}
	or
	\begin{align}
		\label{suppli22}
		u^\ast>0,\ v^\ast=0 \ \text{in}\ [a,b]
	\end{align}
	If \eqref{suppli11} is true, then a solution, $v^\ast$ is exists and is positive of
	\begin{align} 
		\begin{cases}
			\vspace{0.2cm}
			& \frac{\mathsf{d}_2}{\mathsf{d}_1}\left[({\mathsf{d}_1}{ v_x^\ast}-{\alpha_1}{v^\ast})_x\right]+rr_1v^\ast{\left[1-\displaystyle\frac{v^\ast}{K_1}\right]}=0,\quad\displaystyle a<x<b,\ t>0 \\
			\vspace{0.2cm}
			&{\mathsf{d}_1}{v_x^\ast}(a,t)-{\alpha_1}v^\ast(a,t)={\mathsf{d}_1}{v_x^\ast}(b,t)-{\alpha_1}v^\ast(b,t)=0,\quad\displaystyle t>0,
		\end{cases}
	\end{align}
From Lemma \ref{lma1}, it is clear that $0< \displaystyle\frac{ v_x^\ast}{v^\ast}<\displaystyle\frac{\alpha_1}{\mathsf{d}_1}$ within $[a,b]$\ if $\alpha_1>0$. Rewrite the equation above so that 
	\begin{align} 
		\label{1stt}
		\begin{cases}
			\vspace{0.2cm}
			&({\mathsf{d}_1}{v_x^\ast}-{\alpha_1}{v^\ast})_x+rr_1v^\ast{\left[1-\displaystyle\frac{v^\ast}{K_1}\right]}=
			\left(1-\frac{\mathsf{d}_2}{\mathsf{d}_1}\right)\left[({\mathsf{d}_1}{v_x^\ast}-{\alpha_1}{v^\ast})_x\right],\quad\displaystyle a<x<b,\ t>0 \\
			\vspace{0.2cm}
			&{\mathsf{d}_1}{v_x^\ast}(a,t)-{\alpha_1}v^\ast(a,t)={\mathsf{d}_1}{v_x^\ast}(b,t)-{\alpha_1}v^\ast(b,t)=0,\quad\displaystyle t>0,
		\end{cases}
	\end{align}
	Because $\tilde{u}^{\alpha_2}=\displaystyle \frac{u_{\alpha_2}^{\ast}}{\| u_{\alpha_2}^{\ast}\|_{L^\infty}}$, for elliptic equation $\tilde{u}^\ast$ satisfies the following yy standard $L_p$ regularity theory \cite{extra} 
	\begin{align} 
		\label{2ndd}
		\begin{cases}
			\vspace{0.2cm}
			& \left[({\mathsf{d}_1}{ \tilde{u}_x^\ast}-{\alpha_1}{\tilde{u}^\ast})_x\right]+rr_1\tilde{u}^\ast{\left[1-\displaystyle\frac{v^\ast}{K_1}\right]}=0,\quad\displaystyle a<x<b,\ t>0 \\
			\vspace{0.2cm}
			&{\mathsf{d}_1}{\tilde{u}_x^\ast}(a,t)-{\alpha_1}\tilde{u}^\ast(a,t)={\mathsf{d}_1}{\tilde{u}_x^\ast}(b,t)-{\alpha_1}\tilde{u}^\ast(b,t)=0,\quad\displaystyle t>0,
		\end{cases}
	\end{align}
	Again by Lemma \ref{lma1}, we have $\displaystyle\frac{ \tilde{u}_x^\ast}{\tilde{u}^\ast}<\displaystyle\frac{\alpha_1}{\mathsf{d}_1}\Rightarrow{\mathsf{d}_1}{ \tilde{u}_x}<{\alpha_1} \tilde{u}^\ast$ in $[a,b]$\ if $\alpha_1>0;$
Integrate over $[a,b]$ after multiplying the equation \eqref{1stt} by $e^{-\frac{\alpha_1}{\mathsf{d}_1}x}\tilde{u}^\ast$, we obtain
	\begin{multline*}
		\displaystyle \int_{a}^{b}{
			[\mathsf{d}_1v_x^\ast-\alpha_1 v^\ast]_x\cdot e^{-\frac{\alpha_1}{\mathsf{d}_1}x}\tilde{u}^\ast+rr_1v^\ast\left[1-\displaystyle\frac{v^\ast}{K_1}\right]e^{-\frac{\alpha_1}{\mathsf{d}_1}x}\tilde{u}^\ast}dx\\
		=\left(1-\frac{\mathsf{d}_2}{\mathsf{d}_1}\right)\displaystyle \int_{a_1}^{b_1}{	[\mathsf{d}_1 v_x^\ast-\alpha_1 v^\ast]_x\cdot e^{-\frac{\alpha_1}{\mathsf{d}_1}x}\tilde{u}^\ast}dx
	\end{multline*}
	\begin{multline}
		\label{e11}
		\Rightarrow -\displaystyle \int_{a}^{b}{[\mathsf{d}_1 \tilde{u}_x^\ast-\alpha_1 \tilde{u}^\ast]\cdot [\mathsf{d}_1v_x^\ast-\alpha_1 v^\ast]e^{-\frac{\alpha_1}{\mathsf{d}_1}x}}dx+
		\displaystyle\int_{a}^{b}{rr_1v^\ast\left[1-\displaystyle\frac{v^\ast}{K_1}\right]e^{-\frac{\alpha_1}{\mathsf{d}_1}x}\tilde{u}^\ast}dx\\
		=-\left(1-\frac{\mathsf{d}_2}{\mathsf{d}_1}\right)\displaystyle \int_{a}^{b}{[\mathsf{d}_1\tilde{u}_x^\ast-\alpha_1 \tilde{u}^\ast]\cdot [\mathsf{d}_1 v_x^\ast-\alpha_1 v^\ast]e^{-\frac{\alpha_1}{\mathsf{d}_1}x}}dx
	\end{multline}
	Integrate over $[a,b]$ after multiplying the equation \eqref{2ndd} by $e^{-\frac{\alpha_1}{\mathsf{d}_1}x}v^\ast$, we obtain
	\begin{equation*}
		\displaystyle \int_{a}^{b}{
			[\mathsf{d}_1  \tilde{u}_x^\ast-\alpha_1 \tilde{u}^\ast]_x\cdot e^{-\frac{\alpha_1}{\mathsf{d}_1}x}v^\ast+rr_1\tilde{u}^\ast\left[1-\displaystyle\frac{v^\ast}{K_1}\right]e^{-\frac{\alpha_1}{\mathsf{d}_1}x}v^\ast}dx=0
	\end{equation*}
	\begin{multline}
		\label{e22}
		\Rightarrow -\displaystyle \int_{a}^{b}{[\mathsf{d}_1\tilde{u}_x^\ast-\alpha_1 \tilde{u}^\ast]\cdot [\mathsf{d}_1 v_x^\ast-\alpha_1 v^\ast]e^{-\frac{\alpha_1}{\mathsf{d}_1}x}}dx+\displaystyle\int_{a}^{b}{rr_1\tilde{u}^\ast\left[1-\displaystyle\frac{v^\ast}{K_1}\right]e^{-\frac{\alpha_1}{\mathsf{d}_1}x}v^\ast}dx=0
	\end{multline}
	Subtracting equation \eqref{e22} from the equation \eqref{e11}, we get
	$$\left(1-\frac{\mathsf{d}_2}{\mathsf{d}_1}\right)\displaystyle \int_{a}^{b}{[\mathsf{d}_1\tilde{u}_x^\ast-\alpha_1 \tilde{u}^\ast]\cdot [\mathsf{d}_1 v_x^\ast-\alpha_1 v^\ast]e^{-\frac{\alpha_1}{\mathsf{d}_1}x}}dx=0$$
However, this isn't feasible. Therefore, our presumption that a coexistence stable state exists is incorrect. The evidence of assertion $2(a)$ is thus provided. We can also establish the statement $2(b)$ in a similar manner. Let $f:\mathbb{R}\rightarrow\mathbb{R}$ be a non negative continuous function and $\mathbb{P}$ be the set of all $f$ in order to verify assertion $2(c)$. Let's define
	$$\mathbb{Q}:=\mathbb{P}\times(-\mathbb{P})$$
	with non empty interior 
	$$\text{Int}\mathbb{Q}=\text{Int}\mathbb{P}\times(-\text{Int}\mathbb{P})$$
Let's use $$\ll_\mathbb{Q},<_\mathbb{Q},\leq_\mathbb{Q}$$ to represent the partial order relations, which are produced by $\mathbb{Q}$ \cite{competitive system2}.The set of equations \eqref{harvest_free} generates the continuous semi flow for every $\alpha_2\in \sigma:=(\omega_1\alpha_1,\alpha_1)$, which is represented by $\mathcal{Y}_t^{\alpha_2}$. Because $\hat{v}$ depends on $\alpha_2$, we shall use $\hat{v}_{\alpha_2}$ for convenience. Describe
	$$\Gamma_{\alpha_2}=\{(u_0,v_0):(0,\hat{v}_{\alpha_2})\leq_\mathbb{Q}(u_0,v_0)\leq_\mathbb{Q}(\hat{u},0),\quad u_0,v_0\not\equiv 0\}$$
	
	$$X=\{\alpha_2\in\Gamma:(\hat{u},0)\ \text{is\ globally\ asymptotically\ stable\ for} (u_0,v_0)\in \Gamma_{\alpha_2}\}$$
	and
	$$Y=\{\alpha_2\in\Gamma:(0,\hat{v}_{\alpha_2})\ \text{is\ globally\ asymptotically\ stable\ for} (u_0,v_0)\in \Gamma_{\alpha_2}\}$$
The intersection of the sets $X$ and $Y$ is empty. $$(\alpha_1-\epsilon_2,\alpha_1)\subset X\ \text{and}\ (\omega_1\alpha_1,\omega_1\alpha_1+\epsilon_1)\subset Y$$ is evident from statements $2(a)$ and $2(b)$. If it is feasible, let us assume that the system \eqref{harvest_free} does not have a co-existence stable state for $\alpha_2\in \sigma:=(\omega_1\alpha_1,\alpha_1)$. This is taken from \cite{competitive system2} that
	\begin{equation}
		\label{open sets1}
		\sigma:=(\omega_1\alpha_1,\alpha_1)=X\cup Y
	\end{equation}
	Suppose $X$ is open. Let $(u_0^1,v_0^1)$ be a point in $\Gamma_{\alpha_2}$ for any $\alpha_2^0$ which satisfies
	$$(0,\hat{v}_{\alpha_2}^0)\leq_\mathbb{Q}(u_0^1,v_0^1)\leq_\mathbb{Q}(\hat{u},0)$$
	Then there exists $\delta>0$ small such that for $|\alpha_2-\alpha_2^0|\delta$,
	$$(0,\hat{v}_{\alpha_2})\leq_\mathbb{Q}(u_0,v_0)\leq_\mathbb{Q}(\hat{u},0)$$
	Now as $t\rightarrow\infty$, $\mathcal{Y}_t^{\alpha_2^0}(u^0,v^0)\rightarrow(\hat{u},0)\gg_\mathbb{Q}(u_0^1,v_0^1)$.
	Therefore, for a large scale oftime, $t_\infty$
	$$\mathcal{Y}_{t_\infty}^{\alpha_2^0}(u^0,v^0)\gg_\mathbb{Q}(u_0^1,v_0^1)$$
	The continuous dependence of $\mathcal{Y}_{t_\infty}^{\alpha_2^0}(u^0,v^0)$ on $\alpha_2$ implies that there exists $0<\delta^\prime<\delta$ such that for all $|\alpha_2-\alpha_2^0|\delta^\prime$,
	$$\mathcal{Y}_{t_\infty}^{\alpha_2^0}(u^0,v^0)\gg_\mathbb{Q}(u_0^1,v_0^1)\gg_\mathbb{Q}(0,\hat{v}_{\alpha_2})$$
	which implies that for any $\alpha_2\in (\alpha_2^0-\delta^\prime,\alpha_2^0+\delta^\prime)$, $\alpha_2\not\in Y$. Hence
	$$(\alpha_2^0-\delta^\prime,\alpha_2^0+\delta^\prime)\subset X.$$
	Our presumption was correct. In the same way, the set $Y$ is open. The open intervals, however, are the union of two open sets, which is contradictory, as we can see from equation \eqref{open sets1}. Statement $2(c)$ is satisfied by this.
\end{proof}

\section{Numerical Illustrations}\label{section_7}
The primary objective of this section is to present a series of numerical examples that visually illustrate the influence of diffusion, advection, and harvesting on the dynamics of the two-species system. These examples are carefully designed to complement and clarify the theoretical findings discussed in the preceding sections. By incorporating various parameter configurations, we aim to provide a comprehensive understanding of how these factors interact to shape population behavior.
The numerical simulations highlight key stability outcomes, including asymptotic, global, and local stability, which are depicted through graphical representations. These visualizations demonstrate the system's convergence toward steady states under different conditions, offering insights into the stability of semi-trivial and coexistence states. Additionally, the results showcase scenarios of coexistence, where both species persist over time, as well as competitive exclusion, where one species drives the other to extinction.

\begin{figure}[H]
	\centering
	\subfloat[]{\includegraphics[scale=.33]{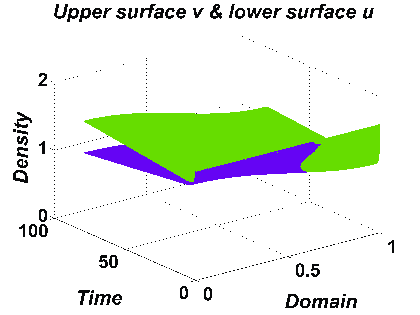}}
	\subfloat[]{\includegraphics[scale=.32]{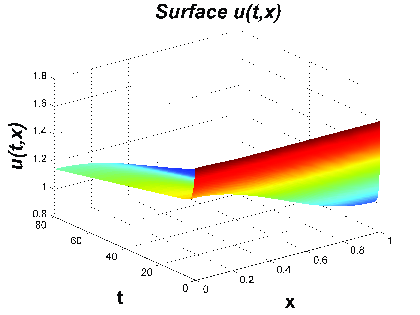}}
	\subfloat[]{\includegraphics[scale=.32]{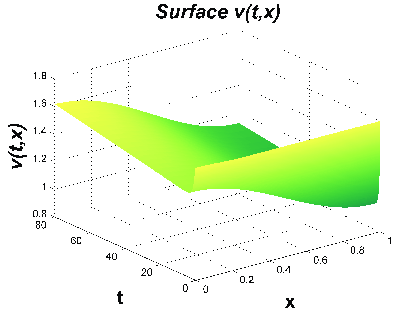}}
	\caption{\small Solution of \eqref{harvest_free} in (a) spatial distribution of $u$ \& $v$, (b) surface distribution of $u$ and (c) surface distribution of $v$, for \ $\mathsf{d}_1=0.08, \mathsf{d}_2=0.07,\alpha_1=0.05,\ \alpha_2=0.04\ \text{and}\ \mu=0.009$\ at time\ $t=2000$.}
	\label{fig:1Dsurf_dist1}
\end{figure}

\begin{figure}[H]
	\centering
	\subfloat[]{\includegraphics[scale=.33]{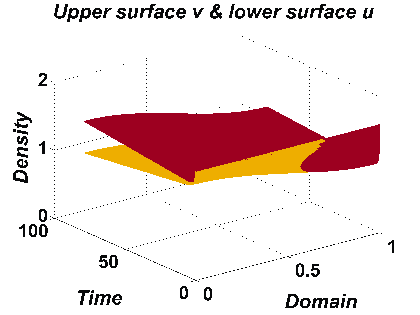}}
	\subfloat[]{\includegraphics[scale=.32]{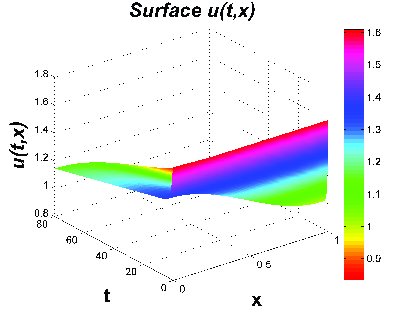}}
	\subfloat[]{\includegraphics[scale=.32]{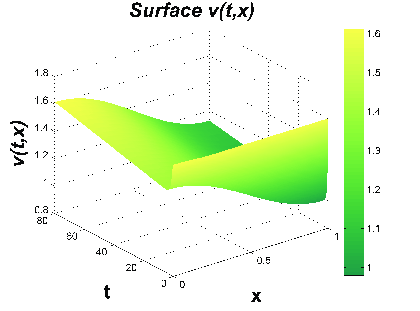}}
	\caption{\small Solution of \eqref{harvest_free} in (a) Spatial distribution of $u$ \& $v$, (b) surface distribution of $u$ and (c) surface distribution of $v$, for \ $\mathsf{d}_1=0.08, \mathsf{d}_2=0.07,\alpha_1=0.05,\ \alpha_2=0.04\ \text{and}\ \mu=0.001$\ at time\ $t=2000$.}
	\label{fig:1Dsurf_dist2}
\end{figure}

Through these simulations, we explore how differences in advection and diffusion rates, combined with harvesting pressures, determine the long-term outcomes of competition. The graphical results effectively bridge the gap between theoretical predictions and practical ecological implications, providing a deeper understanding of the mechanisms governing species persistence, dominance, or extinction. These examples underscore the critical role of parameter variations in determining ecological balance and help validate the robustness of the analytical conclusions drawn in earlier sections.
\begin{Ex}
Consider the system \eqref{harvest_free}, where the carrying capacity is given by 
$K=2.0+\cos(\pi x)$, reflecting spatial variability influenced by periodic environmental heterogeneity, and the harvesting coefficient 
$\mu=0.009\in[0,1)$, which varies over both time and space. These parameter choices allow us to model realistic ecological scenarios where resources and external pressures fluctuate across the habitat.

According to Theorem \ref{th_harvest}, the competition outcome in such a system is fundamentally influenced by the interplay between the diffusion and advection rates of the competing species. Specifically, if a species possesses both higher diffusion and advection rates, the decisive factor becomes the ratio of its advection rate to diffusion rate. The theorem indicates that when the harvesting rate remains within the interval 
$[0,1)$, the species with lower individual advection and diffusion rates will consistently prevail, provided that its advection-to-diffusion ratio is smaller compared to its competitor.


\end{Ex}
The spatial distribution of species 
$u$ and $v$ is depicted in Figures \ref{fig:1Dsurf_dist1}(a) and \ref{fig:1Dsurf_dist2}(a), illustrating a scenario where both species coexist. This coexistence occurs because the diffusion and advection rates of the first species are higher than those of the second species. Additionally, the advection-to-diffusion ratio for the first species is greater than the corresponding ratio for the second species. According to Theorem \ref{th_harvest}, these conditions ensure that both species can persist and share the habitat over time.

The coexistence is further supported by the spatial separation of their distributions, as shown in Figures \ref{fig:1Dsurf_dist1}(b,c) and \ref{fig:1Dsurf_dist2}(b,c), where the individual spatial patterns of species 
$u$ and $v$ are presented separately. These figures highlight how the movement dynamics and interactions between the species, influenced by their respective diffusion and advection rates, allow for the stable partitioning of resources and space within the environment.
This outcome emphasizes the critical role of diffusion and advection rates, as well as their ratios, in determining coexistence. The ability of each species to establish and maintain a presence despite competition depends on these factors, demonstrating the delicate balance of ecological interactions. The spatial distributions depicted in these figures not only validate the theoretical predictions but also provide valuable insights into how movement behaviors and environmental heterogeneity influence the coexistence of competing species in advective systems.

\begin{figure}[H]
	\centering
	\subfloat[]{\includegraphics[scale=.32]{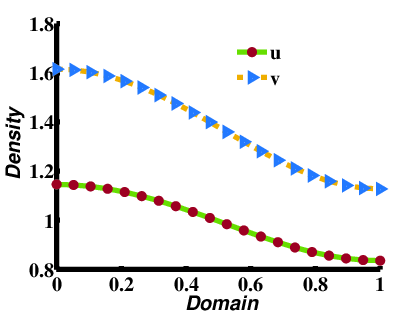}}
	\subfloat[]{\includegraphics[scale=.32]{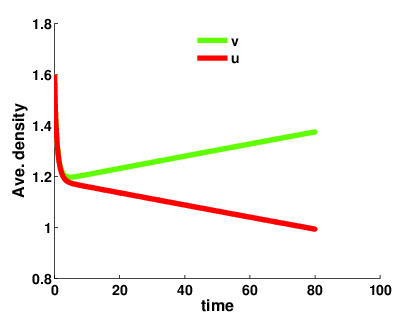}}
	\caption{\small (a) Solution  of \eqref{harvest_free} at time\ $t=2000$ and (b) average solution  of \eqref{harvest_free} at time\ $t=80$, for \ $\mathsf{d}_1=0.08, \mathsf{d}_2=0.07,\alpha_1=0.05,\ \alpha_2=0.04\ \text{and}\ \mu=0.009$\ at time\ $t=2000$.}
	\label{fig:1Dline solution1}
\end{figure}
\begin{figure}[H]
	\centering
	\subfloat[]{\includegraphics[scale=.32]{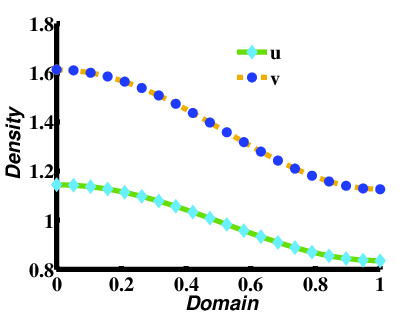}}
	\subfloat[]{\includegraphics[scale=.32]{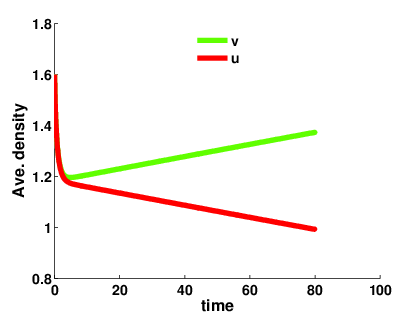}}
	\caption{\small (a) Solution  of \eqref{harvest_free} at time\ $t=2000$ and (b) average solution  of \eqref{harvest_free} at time\ $t=80$, for \ $\mathsf{d}_1=0.08, \mathsf{d}_2=0.07,\alpha_1=0.05,\ \alpha_2=0.04\ \text{and}\ \mu=0.001$\ at time\ $t=80$.}
	\label{fig:1Dline solution2}
\end{figure}

Theorem \ref{th_harvest} predicts that the first species will eventually go extinct due to its higher diffusion and advection rates compared to the second species, as well as a disproportionately larger ratio of these rates. This outcome is visually supported by Figures \ref{fig:1Dline solution1}(a, b) and \ref{fig:1Dline solution2}(a, b), which depict the first species' decline over time. However, Figure \ref{fig:1Dsurf_dist2} illustrates a scenario where both species coexist, as evidenced by their spatial distributions.

Figure \ref{fig:1Dline solution2} illustrates the dynamic relationship between two species over time. While coexistence is apparent on shorter time scales, the first species ultimately faces extinction as time progresses. This transition highlights the temporal nature of species interactions, where short-term persistence does not necessarily guarantee long-term survival.
Initially, both species manage to coexist, suggesting that their population densities are sufficient to sustain each other in the shared habitat.
The shorter time scale likely reflects a period where the first species can counteract losses (e.g., due to diffusion or harvesting) through growth or reproduction.
Over a longer time scale, the first species becomes extinct, indicating that it cannot maintain its population against the pressures it faces.

\begin{Ex}
Consider the system \eqref{harvest_free}, where the carrying capacity is represented by 
$K=2.0+\cos(\pi x)\cos(\pi y)$, introducing spatial heterogeneity due to its dependence on both 
$x$ and $y$. This variability reflects realistic environmental conditions where resource availability changes across different spatial regions. Additionally, the harvesting coefficient 
$\mu=0.03\in[0,1)$ varies over both space and time, accounting for external pressures such as human exploitation or natural predation that fluctuate dynamically.

According to Theorem \ref{th_harvest}, the outcome of competition between two species in this system is primarily influenced by the interplay of their diffusion and advection rates. If one species has both a higher diffusion rate and a higher advection rate compared to its competitor, the critical determinant becomes the ratio of its advection rate to its diffusion rate. The theorem further suggests that when this ratio is smaller for a species, it gains a competitive advantage. Consequently, the species with a lower advection rate and diffusion rate will consistently prevail over its competitor if its advection-to-diffusion ratio is sufficiently reduced.

This result highlights the importance of balancing movement dynamics in heterogeneous environments. A lower advection-to-diffusion ratio allows a species to effectively utilize available resources without being excessively influenced by environmental flows, ensuring stability in its spatial distribution. Conversely, a higher ratio may hinder a species' ability to maintain its presence, as advection-driven movement may cause it to leave favorable areas more quickly than diffusion can compensate. These insights underline the intricate relationship between movement parameters and harvesting effects, offering a nuanced understanding of species competition in spatially and temporally variable environments.	
	
\end{Ex}
\begin{figure}[H]
	\centering
	\subfloat[]{\includegraphics[scale=.33]{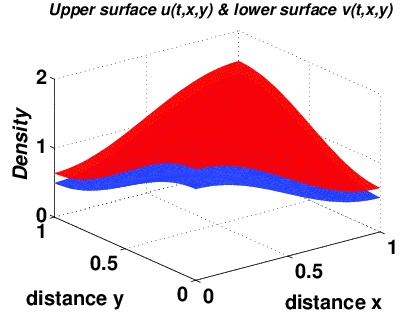}}
	\subfloat[]{\includegraphics[scale=.32]{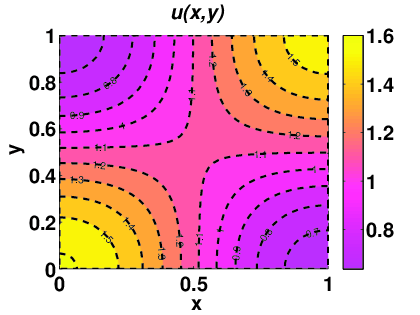}}
	\subfloat[]{\includegraphics[scale=.32]{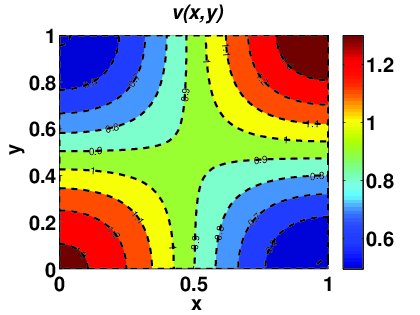}}
	\caption{\small (a) Surface distribution of $u$ \& $v$ of \eqref{harvest_free}, (b) contour plot of $u$ and (c) contour plot of $v$, for \ $\mathsf{d}_1=0.005, \mathsf{d}_2=0.002,\alpha_1=0.002,\ \alpha_2=0.0018\ \text{and}\ \mu=0.03$\ at time\ $t=80$.}
	\label{fig:2Dsurf_dist1}
\end{figure}
As shown in Figure \ref{fig:2Dsurf_dist1}(a), the dynamics of the system reveal that the two species can coexist when the diffusion and advection rates of the second species are smaller than those of the first species. This coexistence is consistent with the results stated in Theorem \ref{th_harvest}(c), which asserts that minor perturbations in the second species' advection rate can facilitate stable coexistence under these conditions. The reduced advection and diffusion rates of the second species allow it to occupy specific niches within the environment without being entirely outcompeted by the more mobile first species.

The contour plots in Figures \ref{fig:2Dsurf_dist1}(b) and \ref{fig:2Dsurf_dist1}(c) provide additional insights into the spatial distribution of the two species. These plots illustrate that the density profile of the first species $(u)$ is consistently higher than that of the second species $(v)$. This indicates that, while the second species is able to coexist with the first, its population density is relatively lower due to its reduced capacity to explore and utilize the habitat effectively compared to the first species, which has higher diffusion and advection rates.

These observations highlight the intricate balance between movement dynamics and competitive interactions in determining coexistence. The ability of the second species to persist despite its disadvantages in diffusion and advection rates underscores the stabilizing role of minor perturbations and localized resource availability. The spatial patterns depicted in the contour plots further emphasize how differences in movement behaviors and environmental interactions shape the density profiles and overall coexistence of the species within a heterogeneous habitat.


\begin{figure}[H]
	\centering
	\subfloat[]{\includegraphics[scale=.33]{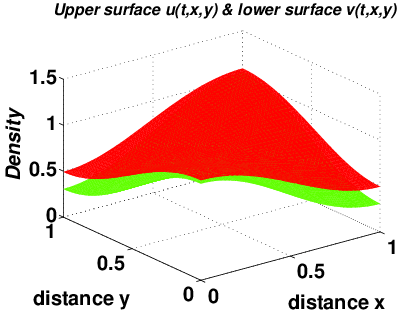}}
	\subfloat[]{\includegraphics[scale=.32]{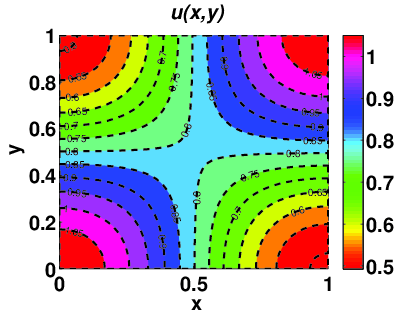}}
	\subfloat[]{\includegraphics[scale=.32]{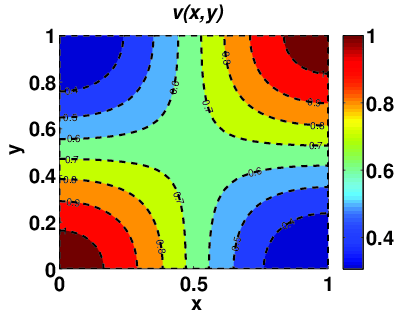}}
	\caption{\small (a) Surface distribution of $u$ \& $v$ of \eqref{harvest_free}, (b) contour plot of $u$ and (c) contour plot of $v$, for \ $\mathsf{d}_1=0.002, \mathsf{d}_2=0.001,\alpha_1=0.001,\ \alpha_2=0.0006\ \text{and}\ \mu=0.3$\ at time\ $t=80$.}
	\label{fig:2Dsurf_dist2}
\end{figure}
According to Theorem \ref{th_harvest}(c), coexistence is seen in Figure \ref{fig:2Dsurf_dist2}(a) if the second species' diffusion and advection rate is smaller and \eqref{c_2} holds compared to the first species, where the harvesting rate is $\mu=0.3$. This is true for small perturbations on the second species' advection rate. Additionally, the contour plots of $u$ and $v$ are displayed in Figures \ref{fig:2Dsurf_dist2}(b,c), and it is evident that the density of two species remains constant due to coexistence.

\begin{figure}[H]
	\centering
	\subfloat[]{\includegraphics[scale=.33]{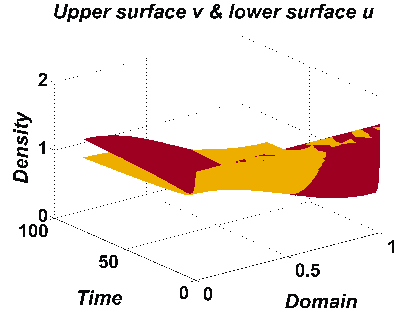}}
	\subfloat[]{\includegraphics[scale=.32]{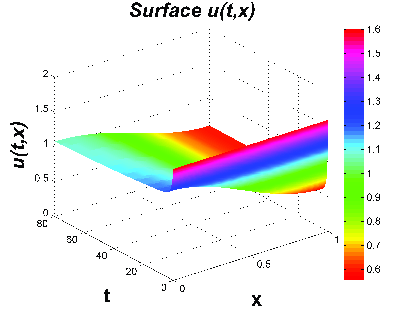}}
	\subfloat[]{\includegraphics[scale=.32]{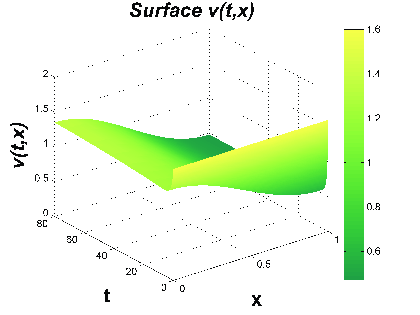}}
	\caption{\small (a) Spatial distribution of $u$ \& $v$, (b) surface solution of $u$ and (c) surface solution of $v$, of \eqref{harvest_free} for \ $\mathsf{d}_1=0.002, \mathsf{d}_2=0.001,\alpha_1=0.001,\ \alpha_2=0.0006\ \text{and}\ \mu=0.3$\ at time\ $t=2000$.}
	\label{fig:1Dsurf_dist3}
\end{figure}
The spatial distribution of species $u$ and $v$ in Figure \ref{fig:1Dsurf_dist3}(a) shows that, in accordance with Theorem \ref{th_harvest}(c), coexistence is feasible. Individual species $u$ and $v$ spatial distributions are shown individually in Figure \ref{fig:1Dsurf_dist3}(b,c). Rivers are dynamic environments influenced by factors such as water pollution, dam construction, and climate change. Reaction-diffusion-advection models, accounting for spatial heterogeneity, can predict how such changes may alter the flow of nutrients, toxins, or temperature gradients in the river, impacting animal health and biodiversity.

\begin{figure}[H]
	\centering
	\subfloat[]{\includegraphics[scale=.32]{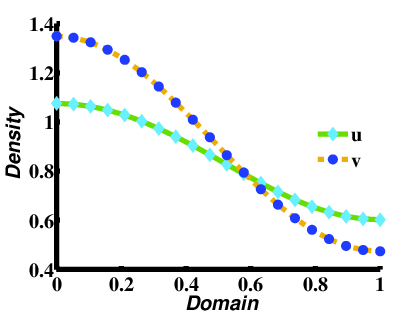}}
	\subfloat[]{\includegraphics[scale=.32]{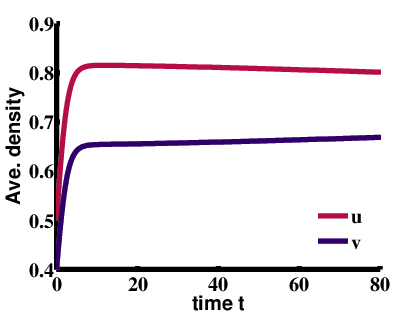}}
	\caption{\small (a) Solution  of \eqref{harvest_free} at time\ $t=2000$ and (b) average solution  of \eqref{harvest_free} at time\ $t=80$, for \ $\mathsf{d}_1=0.002, \mathsf{d}_2=0.001,\alpha_1=0.001,\ \alpha_2=0.0006\ \text{and}\ \mu=0.3$\ at time\ $t=80$.}
	\label{fig:1Dline solution3}
\end{figure}
Figures \ref{fig:1Dline solution3}(a) and \ref{fig:1Dline solution3}(b) demonstrate that when the diffusion and advection rates of the first species are higher, and the condition given by \eqref{c_2} is satisfied, coexistence between the two species can occur. This coexistence is robust to small perturbations in the advection rate of the second species, suggesting a dynamic equilibrium that allows both species to persist under these specific conditions. The higher diffusion and advection rates of the first species facilitate its spread and ability to occupy diverse regions of the domain. However, these advantages must be balanced to prevent excessive dispersal that could otherwise dilute its population.
The condition \eqref{c_2} acts as a critical threshold that defines the parameter space where coexistence is possible. These results can inform strategies to enhance coexistence in multi-species systems, such as designing habitats that optimize dispersal conditions for vulnerable species.
By ensuring environmental factors align with thresholds like \eqref{c_2}, managers can create conditions conducive to long-term biodiversity.

\begin{figure}[H]
	\centering
	\subfloat[]{\includegraphics[scale=.33]{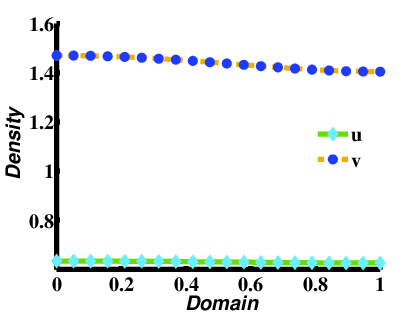}}
	\subfloat[]{\includegraphics[scale=.32]{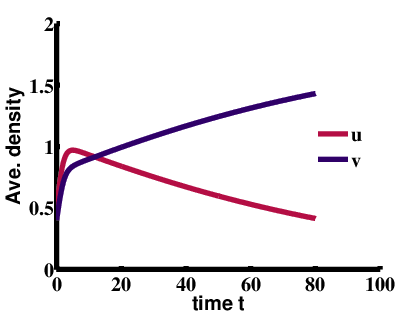}}
	\caption{\small (a) Solution  of \eqref{harvest_free} at time\ $t=2000$ and (b) average solution  of \eqref{harvest_free} at time\ $t=80$, for \ $\mathsf{d}_1=3, \mathsf{d}_2=0.8,\alpha_1=0.7,\ \alpha_2=0.03\ \text{and}\ \mu=0.1$\ at time\ $t=80$.}
	\label{fig:1Dline solution4}
\end{figure}
Figures \ref{fig:1Dline solution4}(a) and \ref{fig:1Dline solution4}(b) depict the dynamics of two species under the influence of harvesting and diffusion, respectively. The first species faces extinction due to the combined effects of high harvesting pressure and an excessively high diffusion rate, as predicted by Theorem \ref{th_harvest}. In contrast, the second species persists, demonstrating its resilience under the given conditions. This observation highlights the importance of understanding species-specific responses to environmental pressures, such as harvesting and dispersal dynamics, which are critical for designing adaptive management strategies aimed at conserving vulnerable species.

The second species benefits from a relatively lower diffusion rate and harvesting pressure, allowing it to maintain sufficient local population densities for survival.
The persistence of the second species underscores its ability to adapt to and withstand the combined effects of these environmental pressures. The dynamics illustrated in the figure are highly relevant for conservation biology and resource management.
By understanding how species respond to different harvesting intensities and diffusion rates, adaptive strategies can be developed to mitigate the risk of extinction for vulnerable species.
For instance, reducing harvesting rates or creating protected areas where diffusion effects are minimized could help stabilize at-risk populations. Theorem \ref{th_harvest} provides a theoretical framework for predicting species outcomes based on specific harvesting and diffusion parameters.
\begin{figure}[H]
	\centering
	\subfloat[]{\includegraphics[scale=.32]{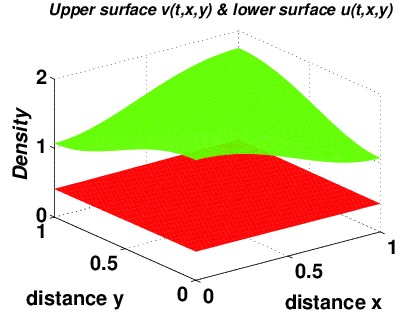}}
	\subfloat[]{\includegraphics[scale=.32]{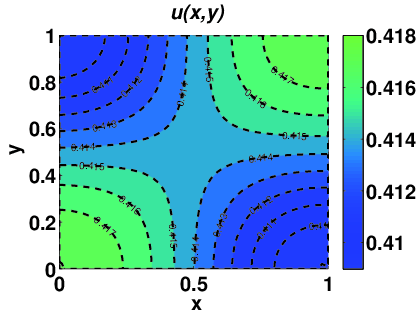}}
	\subfloat[]{\includegraphics[scale=.32]{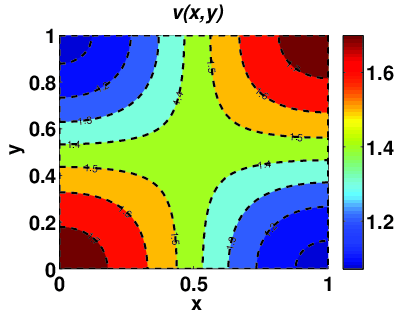}}	
	\caption{\small (a) Surface distribution of $u$ \& $v$ of \eqref{harvest_free}, (b) contour plot of $u$ and (c) contour plot of $v$, for \ $\mathsf{d}_1=3, \mathsf{d}_2=0.8,\alpha_1=0.7,\ \alpha_2=0.03\ \text{and}\ \mu=0.1$\ at time\ $t=80$.}
	\label{fig:2Dsurf_dist3}
\end{figure}
Figure \ref{fig:2Dsurf_dist3}(a) illustrates the surface distribution of the model described in \eqref{harvest_free}. Here, the first species ultimately faces extinction due to its significantly higher diffusion rate compared to the second species. This disparity in diffusion rates causes the first species to spread out too rapidly, reducing its local density to unsustainable levels. Figures \ref{fig:2Dsurf_dist3}(b) and \ref{fig:2Dsurf_dist3}(c) present contour plots for the two species, $u$ and $v$, respectively, highlighting the dominance of the second species. The density profile of the second species $(v)$ is consistently higher than that of the first species $(u)$. Interestingly, despite the differing densities, the spatial distributions suggest that both species are gathered in similar regions, indicating overlapping habitat use.

The first species $(u)$ suffers from an extremely high diffusion rate, which leads to rapid dispersal across the domain.
While higher diffusion rates might initially help the species colonize new areas, the inability to maintain concentrated populations in any one region leads to unsustainably low local densities (Figure \ref{fig:2Dsurf_dist3}(a)).
Over time, this dynamic causes the population of the first species to dwindle to extinction.

The contour plots for the two species show stark differences in their density profiles.
The second species $(v)$, with its lower diffusion rate, maintains higher local densities, allowing it to dominate the system.
The higher density of the second species is a direct consequence of its ability to retain populations in concentrated areas, minimizing losses and outcompeting the first species (Figures \ref{fig:2Dsurf_dist3}(b, c)).

\begin{figure}[H]
	\centering
	\subfloat[]{\includegraphics[scale=.32]{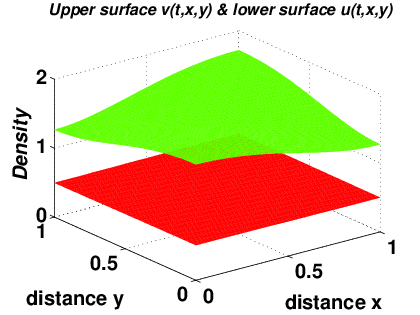}}
	\subfloat[]{\includegraphics[scale=.32]{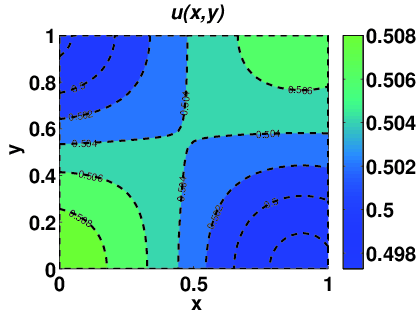}}
	\subfloat[]{\includegraphics[scale=.1]{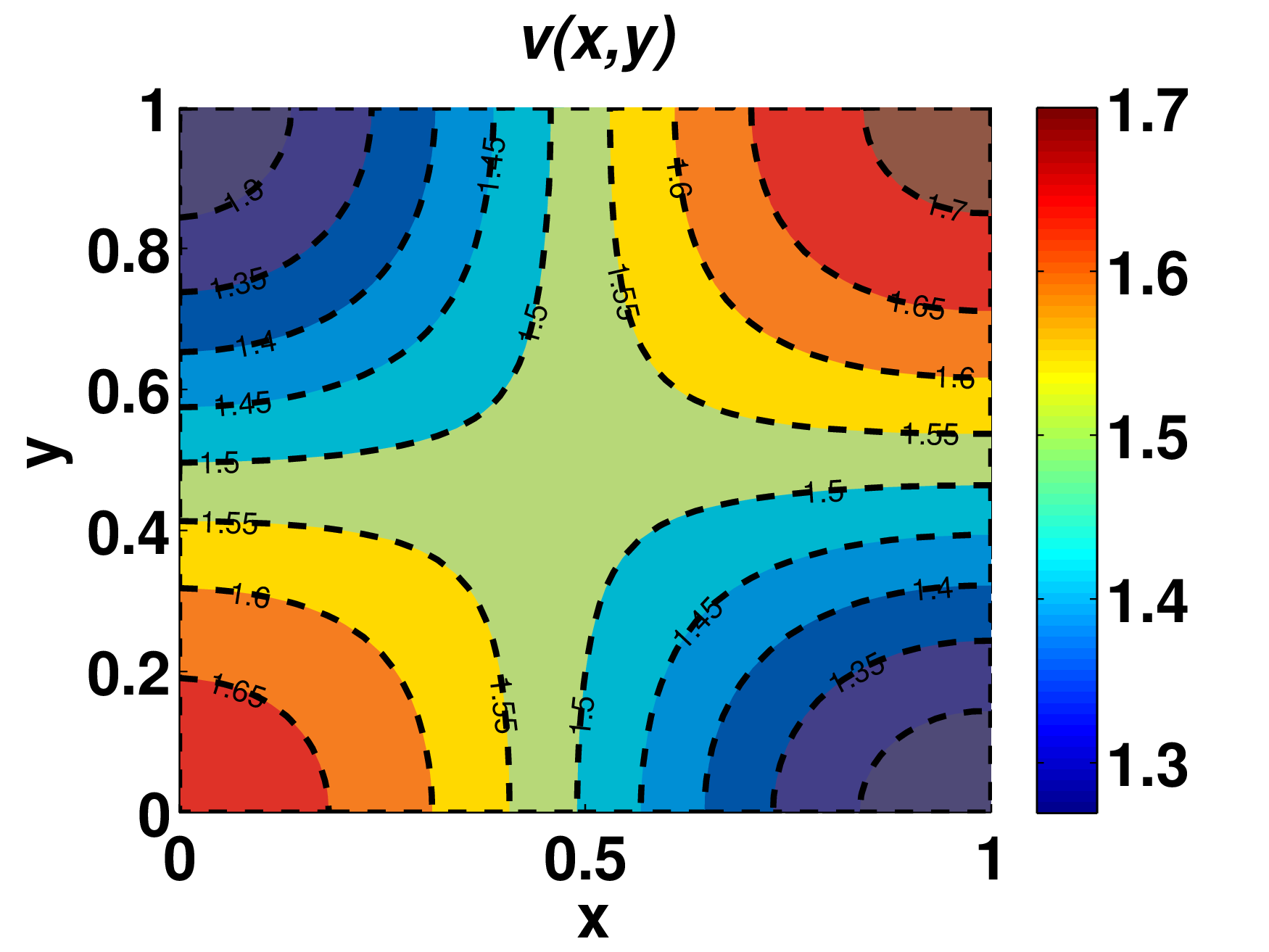}}	
	\caption{\small (a) Surface distribution of $u$ \& $v$ of \eqref{harvest_system}, (b) contour plot of $u$ and (c) contour plot of $v$, for \ $\mathsf{d}_1=1.0, \mathsf{d}_2=1.0,\alpha_1=1.0,\ \alpha_2=1.0,\ \mu_1=0.01\ \text{and}\ \mu_2=0.0076$\ at time\ $t=80$.}
	\label{fig:2Dsurf_dist4}
\end{figure}
Figure \ref{fig:2Dsurf_dist4}(a) illustrates the spatial distribution of two species when the harvesting rate for the first species exceeds that of the second species, leading to coexistence between the two populations. In contrast, Figure \ref{fig:2Dsurf_dist4}(b) and \ref{fig:2Dsurf_dist4}(c) present contour plots where the density of the second species surpasses that of the first. This outcome is observed under varying harvesting rates. For this analysis, we assume equal diffusion and advection rates for both species. The figures are generated for a relatively shorter time scale, specifically at $t=80$. Over longer time periods, the dynamics shift, resulting in the extinction of the first species while the second species continues to persist.

\begin{itemize}
	\item Coexistence Dynamics (Figure \ref{fig:2Dsurf_dist4}(a)):
	The first species, despite being subjected to higher harvesting rates, still maintains a presence across the surface alongside the second species.
	This coexistence may reflect a balance between the growth, dispersal, and mortality processes that allow the first species to withstand its higher harvesting pressure.
	The spatial distribution likely highlights the regions where the two species overlap, creating a scenario of shared habitat.
	\item Dominance of the Second Species (Figures \ref{fig:2Dsurf_dist4}(b, c)):
	The second species, benefiting from lower harvesting rates, achieves a higher population density than the first species.
	The contour plots may showcase how the second species takes advantage of favorable conditions, potentially outcompeting the first in certain areas.
	This dynamic arises solely from differences in harvesting pressures, assuming other factors like dispersal (diffusion and advection) remain equal.
	\item Role of Diffusion and Advection:
	By assuming equal rates for diffusion and advection, the spatial spread of both species is determined purely by their harvesting rates and initial distributions.
	This simplifies the system, isolating the effects of harvesting rates on population dynamics.
	\item Temporal Dynamics:
	At $t=80$, the figures capture a snapshot of the species’ distribution, representing an intermediate stage before longer-term outcomes become apparent.
	Over extended periods, the disparity in harvesting rates leads to the eventual extinction of the first species. The second species, with its comparatively lower harvesting pressure, establishes dominance and persists in the habitat.
\end{itemize}
This analysis underscores the critical role that harvesting rates play in determining species distribution, density, and long-term viability, even when other ecological parameters such as diffusion and advection are held constant.

\section{Conclusion}
The dynamics of reaction-diffusion-advection models provide valuable insights into the interactions between species, including predator-prey relationships, competition for limited resources, and symbiotic partnerships \cite{k1}. These models are particularly significant in understanding how spatial heterogeneity within ecosystems, such as river systems, influences these interactions. Such insights are crucial for predicting species survival and behavior, especially under changing environmental conditions that can disrupt ecological balance.
In this study, we analyzed a reaction-diffusion-advection model incorporating harvesting efforts for two competing species in a heterogeneous advective environment. The model operates under zero Neumann boundary conditions, ensuring that individuals cannot cross the upstream and downstream boundaries, reflecting realistic constraints in river systems. For our analysis, we assumed that the harvesting rate lies within the interval 
$[0,1)$. The primary focus was to examine the global stability of the two competing species under these conditions.

We first established the model's foundational properties, including the existence, uniqueness, and positivity of its solutions. These ensure that the model behaves in a biologically realistic manner, with population densities remaining non-negative and well-defined. Furthermore, we demonstrated the local stability of two semi-trivial steady states, representing scenarios where one species dominates while the other is excluded. We then explored the non-existence of a coexistence steady state under specific non-trivial assumptions, highlighting the conditions under which one species inevitably outcompetes the other.
The central result of this study is encapsulated in Theorem \ref{th_harvest}, developed using the concept of monotone dynamical systems and two key hypotheses. The theorem establishes that the competition's outcome is determined by the ratio of advection rate to diffusion rate, particularly when both species employ dispersion strategies characterized by low harvesting, diffusion, and advection rates. Specifically, the species with the lower diffusion and advection rates will always prevail if its advection-to-diffusion ratio is reduced. This finding underscores the critical role of movement dynamics and environmental interactions in shaping competitive outcomes.

Despite the first species being on the path to extinction, the contour plots reveal a significant degree of overlap in the regions where the two species are found (Figure \ref{fig:2Dsurf_dist3}).
This suggests that both species are drawn to similar environmental conditions or resources, leading to a shared habitat despite the disparities in their diffusion rates.
The gathering of both species in similar areas may be an artifact of the model’s assumptions, such as shared environmental gradients or mutual attraction to favorable conditions. Diffusion rates play a crucial role in shaping species distributions, especially in competitive scenarios.
A high diffusion rate, as seen in the first species, can be detrimental when it prevents populations from remaining concentrated enough to sustain reproduction and survival.
Conversely, the second species benefits from its slower diffusion rate, which allows it to establish and maintain stable population densities. This analysis underscores the interplay between diffusion rates and species interactions, demonstrating how these factors influence spatial distribution, density profiles, and long-term survival outcomes in ecological systems. Applying such theoretical insights to real-world scenarios can inform policy decisions and optimize management practices. The results emphasize the delicate balance between species dispersal and local population dynamics.
High diffusion rates, while seemingly advantageous, can have counterintuitive and detrimental effects, especially when combined with external pressures like harvesting.
Adaptive management approaches should consider these interactions to ensure the long-term sustainability of ecosystems.

Beyond its theoretical contributions, this study has practical implications for riverine ecology. The findings offer a quantitative framework for understanding how physical and biological processes interact in spatially complex river environments. By modeling the effects of spatial heterogeneity, such as variations in flow rates, water temperature, or habitat structure, researchers can identify critical zones for conservation efforts. These models can guide interventions like habitat restoration, pollution control, or adaptive harvesting strategies to maximize their impact on the survival and sustainability of river ecosystems.
Ultimately, this work supports the conservation of riverine species by providing actionable insights into the interplay of environmental and biological factors, helping ensure the resilience and sustainability of these vital ecosystems.

\section*{Acknowledgements}
The research by M. Kamrujjaman was partially supported by the
University Grants Commission (UGC), and the University of Dhaka,
Bangladesh.

\section*{Declaration of competing interest}
The authors declare no competing interests. 

\section*{Data availability}
No data used in this study. 

\section*{Ethical approval}
N/A. 

\section*{CRediT authorship contribution statement} 
Mayesha Sharmim Tisha: Conceptualization, Data curation, Formal analysis, Methodology, Software, Validation, Writing– original draft.\\
Md. Kamrujjaman: Conceptualization, Validation, Resources, Investigation,  Funding acquisition, Software, Supervision, Writing– original draft, Writing– review \& editing.\\
All authors have read and agreed to the published version of the manuscript.


\end{document}